\documentclass[draft]{amsart}
\usepackage[all]{xy}
\usepackage{enumerate}
\usepackage{amssymb}
\usepackage{color}

\newcommand{\ie}{{\em i.e.\ }}
\newcommand{\eg}{{\em e.g.\ }}
\newcommand{\cf}{{\em cf.\ }}

\newtheorem{theorem}{Theorem}
\newtheorem*{thm}{Theorem}
\newtheorem{lemma}{Lemma}
\newtheorem{proposition}{Proposition}
\newtheorem{corollary}{Corollary}

\theoremstyle{definition}
\newtheorem{definition}{Definition}

\newtheorem{example}{Example}
\newtheorem{remark}{\textbf{Remark}}

\numberwithin{equation}{section}

\newcommand{\N}{\mathbf{N}}
\newcommand{\Z}{\mathbf{Z}}

\newcommand{\ko}{\: , \;}

\newcommand{\we}{\wedge}

\renewcommand{\tilde}[1]{\widetilde{#1}}

\newcommand{\ra}{\rightarrow}

\newcommand{\arr}[1]{\stackrel{#1}{\rightarrow}}

\newcommand{\opname}[1]{\operatorname{\mathsf{#1}}}

\newcommand{\Mod}{\opname{Mod}\nolimits}
\newcommand{\Sum}{{\mbox{Sum}}}

\newcommand{\Tria}{{\mbox{Tria}}}
\newcommand{\Susp}{{\mbox{Susp}}}
\newcommand{\aisle}{{\mbox{aisle}}}

\newcommand{\id}{\mathbf{1}}
\newcommand{\R}{\mathbf{R}}
\renewcommand{\L}{\mathbf{L}}

\newcommand{\colim}{\colim}
\newcommand{\lid}{\varinjlim}

\newcommand{\Mcolim}{\opname{Mcolim}}
\newcommand{\can}{\opname{can}}

\newcommand{\op}[1]{\opname{#1}\nolimits}
\newcommand{\Zy}[1]{\op{Z}^{#1}}

\renewcommand{\H}[1]{{H}^{#1}}

\newcommand{\ca}{{\mathcal A}}
\newcommand{\cb}{{\mathcal B}}
\newcommand{\cc}{{\mathcal C}}
\newcommand{\cd}{{\mathcal D}}

\newcommand{\ch}{{\mathcal H}}

\newcommand{\cl}{{\mathcal L}}
\newcommand{\cm}{{\mathcal M}}

\newcommand{\cp}{{\mathcal P}}

\newcommand{\cq}{{\mathcal Q}}

\newcommand{\cu}{{\mathcal U}}
\newcommand{\cv}{{\mathcal V}}

\newcommand{\cx}{{\mathcal X}}
\newcommand{\cy}{{\mathcal Y}}
\newcommand{\cz}{{\mathcal Z}}
\renewcommand{\phi}{\varphi}

\newcommand{\Hom}{\opname{Hom}}
\newcommand{\End}{\opname{End}}

\newcommand{\Ext}{\opname{Ext}}

\newcommand{\cone}{\opname{Cone}\nolimits}

\begin{document}
\title{Lifting and restricting recollement data}

\author{Pedro Nicol\'{a}s}
\address{Departamento de Matem\'{a}ticas, Universidad de Murcia, Aptdo. 4021, 30100, Espinardo, Murcia, Espa\~na}

\email{pedronz@um.es}
\thanks{The authors have been partially supported by research projects from the D.G.I. of the Spanish Ministry of Education and theFundaci\'{o}n S\'{e}neca of Murcia, with a part of FEDER funds. The first author has been also supported by the MECD grant AP2003-2896.}
\author{Manuel Saor\'{i}n}
\address{Departamento de Matem\'{a}ticas, Universidad de Murcia, Aptdo. 4021, 30100, Espinardo, Murcia, Espa\~na}\email{msaorinc@um.es}

\date{March 31, 2008}



\begin{abstract}
We study the problem of lifting and restricting TTF triples (equivalently, recollement data) for a certain wide type of triangulated categories. This, together with the parametrizations of TTF triples given in \cite{NicolasSaorin2007c}, allows us to show that many well-known recollements of right bounded derived categories of algebras are restrictions of recollements in the unbounded level, and leads to criteria to detect recollements of general right bounded derived categories. In particular, we give in Theorem \ref{parametrization right bounded recollements} necessary and sufficient conditions for a \emph{right bounded} derived category of a differential graded(=dg) category to be a recollement of right bounded derived categories of dg categories. In Theorem \ref{ourKonig} we consider the particular case in which those dg categories are just ordinary algebras.\end{abstract}

\maketitle
\tableofcontents

\section{Introduction}

\subsection{Motivations}

\emph{Torsion torsion-free(=TTF) triples} are important in the
theory of abelian categories (in particular, categories of
modules), \cf for instance \cite{Stenstrom}. It turns out that TTF
triples still `make sense' in the theory of triangulated
categories and that they are also important for they are in
bijection with \emph{recollement data} (\cf subsection \ref{TTF
triples and recollement data}) and, in many cases, with
\emph{smashing subcategories} (\cf \cite[Proposition
4.4.14]{NicolasTesis}, \cite[Corollary 2.4]{NicolasSaorin2007c}).

Once the problem of parametrizing TTF triples on perfectly
generated triangulated categories (in particular, unbounded
derived categories of small dg categories) has been essentially
solved in \cite{NicolasSaorin2007c}, we study here the problem of
lifting and restricting TTF triples for certain natural full
triangulated subcategories which generalize the subcategory of the
derived category of an algebra formed by the complexes with right
bounded cohomology. A byproduct of our results is  an `unbounded'
approach to S.~K\"{o}nig's work \cite{Konig1991}.

\subsection{Outline of the paper}

In section \ref{Notation and preliminary results}, we fix some terminology and recall some results on triangulated categories. Also, we introduce the \emph{right bounded derived category} of a small dg category. In section \ref{Lifting of TTF triples}, we study the problem of \emph{lifting} a TTF triple from a certain full triangulated subcategory $\cd'$ of a triangulated category $\cd$ with small coproducts and a set of generators contained in $\cd'$. In subsection \ref{General criterion}, we consider the general case, and in subsection \ref{`Right bounded' triangulated subcategories} we focus on the case in which $\cd'$ is a kind of `right bounded' triangulated subcategory of $\cd$. In section \ref{Restriction of TTF triples}, we study the problem of \emph{restricting} TTF triples. The general criterion (\cf subsection \ref{General criterion bis}) was already given by A. A. Beilinson, J. Bernstein and P. Deligne in their seminal paper \cite{BeilinsonBernsteinDeligne}. In subsection \ref{`Right bounded' triangulated subcategories bis}, we deduce the criterion for the case of a `right bounded' triangulated subcategory. This allows us to regard, in Example \ref{restriction in examples}, some well-known recollements of right bounded derived categories of algebras as restrictions of a recollement induced at the unbounded level by a \emph{homological epimorphism} of the form $A\ra A/I$ where $I$ is a two-sided ideal of the algebra $A$. With the help of the former sections, we study in section \ref{Recollements of right bounded derived categories} the problem of giving necessary and sufficient conditions for a right bounded derived category of a dg category to be a recollement of right bounded derived categories of dg categories. This leads us to inspect in subsection \ref{bounds} some `boundness' conditions for sets of objects of a right bounded derived category of a dg category. In subsection \ref{Recollement of general right bounded derived categories}, we first give a general criterion (\cf Theorem \ref{parametrization right bounded recollements}) and then a criterion (\cf Corollary \ref{parametrization right bounded hn recollements}) for the case when the `glued' dg categories have cohomology concentrated in non-positive degrees. This allows us to deduce, in subsection \ref{Recollement of right bounded derived categories of algebras}, a set of necessary and sufficient conditions for the right bounded derived category of an ordinary algebra to be a recollement of right bounded derived categories of ordinary algebras. A result in that direction already appeared in S.~K\"{o}nig's paper \cite[Theorem 1]{Konig1991}, but we show in section \ref{More than an exceptional object} that stronger assumptions are needed in order S.~K\"{o}nig's theorem to be true in general.
\section{Notation and preliminary results}\label{Notation and preliminary results}
\subsection{Notation}

Unless otherwise stated, $k$ will be a commutative (associative, unital) ring and every additive category will be assumed to be $k$-linear. We will only work with unital algebras and unital modules. We denote by $\Mod k$ the category of $k$-modules. Given a class $\cq$ of objects of an additive category $\cd$, we denote by $\cq^{\bot_{\cd}}$, or $\cq^{\bot}$ if the category $\cd$ is clear, the full subcategory of $\cd$ formed by the objects $M$ which are \emph{right orthogonal} to every object of $\cq$, \ie such that $\cd(Q,M)=0$ for all $Q$ in $\cq$. Dually for $^{\bot_{\cd}}\cq$. When $\cd$ is a triangulated category, the \emph{shift functor} will be denoted by $?[1]$, and its quasi-inverse will be denoted by $?[-1]$. When we use expression like ``all the shifts'' or ``closed under shifts'' and so on, we will mean ``shifts in both directions'', that is to say, we will refer to the $n$th power $?[n]$ of $?[1]$ for all the integers $n\in\Z$. In case we want to consider another situation (\eg non-negative shifts $?[n]\ko n\geq 0$) this will be said explicitly.If $\cq$ is a class of objects of a triangulated category $\cd$:
\begin{enumerate}
\item $\cq^+$ will be the class of all non-negative shifts of objects of $\cq$.\item $\Sum_{\cd}(\cq)$, or $\Sum(\cq)$ if $\cd$ is clear, will be the class of all small coproducts of objects of $\cq$.\item $\aisle_{\cd}(\cq)$, or $\aisle(\cq)$ if $\cd$ is clear, will be the smallest \emph{aisle} (\cf \cite[Definition 1.1]{KellerVossieck88b}) in $\cd$ containing $\cq$. Notice that $\aisle_{\cd}(\cq)$ might not exist since the intersection of aisles might not be an aisle, but if it does then it is closed under small coproducts.\item $\Susp_{\cd}(\cq)$, or $\Susp(\cq)$ if $\cd$ is clear, will be the smallest full \emph{suspended} subcategory (\cf \cite[subsection 1.1]{KellerVossieck1987}) of $\cd$ containing $\cq$ and closed under small coproducts.\item $\Tria_{\cd}(\cq)$, or $\Tria(\cq)$ if $\cd$ is clear, will be the smallest full triangulated subcategory of $\cd$ containing $\cq$ and closed under small coproducts.\end{enumerate}
If $\cu$ and $\cv$ are two classes of objects of a triangulated category $\cd$, then $\cu*\cv$ is the class of \emph{extensions} of objects of $\cv$ by objects of $\cu$, \ie the class formed by those objects $M$ occuring in a triangle\[U\ra M\ra V\ra U[1]
\]
of $\cd$ with $U\in\cu$ and $V\in\cv$. Notice that the operation $*$ is associative. For each natural number $n\geq 0$ the objects of $\cu^{*n}:=\cu*\overset{n\text{ times}}{\dots}*\cu$ are called \emph{$n$-fold extensions of length $n$} of objects of $\cu$.
We will use without explicit mention the bijection between t-structures on a triangulated category $\cd$ and aisles in $\cd$, proved by B.~Keller and D.~Vossieck in \cite{KellerVossieck88b}. If $(\cu,\cv[1])$ is a t-structure on a triangulated category $\cd$, we denote by $u:\cu\hookrightarrow\cd$ and $v:\cv\hookrightarrow\cd$ the inclusion functors, by $\tau_{\cu}$ a right adjoint to $u$ and by $\tau^{\cv}$ a left adjoint to $v$.
\subsection{TTF triples and recollement data}\label{TTF triples and recollement data}
A \emph{torsion torsionfree(=TTF) triple} on a triangulated category $\cd$ is a triple $(\cx,\cy,\cz)$ of full subcategories of $\cd$ such that $(\cx,\cy)$ and $(\cy,\cz)$ are t-structures on $\cd$. Notice that, in particular, $\cx\ko \cy$ and $\cz$ are full triangulated subcategories of $\cd$. It is well known that TTF triples are in bijection with (suitable equivalence classes of) \emph{recollement data} (\cf \cite[1.4.4]{BeilinsonBernsteinDeligne}, \cite[subsection 9.2]{Neeman2001}, \cite[subsection 4.2]{NicolasTesis}). For the convenience of the reader we recall how this bijection works. If\[\xymatrix{ \cd_{F}\ar[r]^{i_{*}} & \cd\ar@/_-1pc/[l]^{i^!}\ar@/_1pc/[l]_{i^*}\ar[r]^{j^*} & \cd_{U}\ar@/_-1pc/[l]^{j_{!}}\ar@/_1pc/[l]_{j_{*}}}
\]
expresses $\cd$ as a recollement of $\cd_{F}$ and $\cd_{U}$, then
\[(j_{!}(\cd_{U}),i_{*}(\cd_{F}),j_{*}(\cd_{U}))
\]
is a TTF triple on $\cd$, where by $j_{!}(\cd_{U})$ we mean the essential image of $j_{!}$, and analogously with the other functors. Conversely, if $(\cx,\cy,\cz)$ is a TTF triple on $\cd$, then $\cd$ is a recollement of $\cy$ and $\cx$ as follows:\[\xymatrix{ \cy\ar[r]^{y} & \cd\ar@/_-1pc/[l]^{\tau_{\cy}}\ar@/_1pc/[l]_{\tau^\cy}\ar[r]^{\tau_{\cx}} & \cx,\ar@/_-1pc/[l]^{x}\ar@/_1pc/[l]_{z\tau^{\cz}x}}
\]
Notice that for a TTF triple $(\cx,\cy,\cz)$ the compositions $\cx\arr{x}\cd\arr{\tau^{\cz}}\cz$ and $\cz\arr{z}\cd\arr{\tau_{\cx}}\cx$ are mutually quasi-inverse triangle equivalences (\cf \cite[Lemma 1.6.7]{NicolasTesis}).
\subsection{(Super)perfectness and compactness}

An object $P$ of a triangulated category $\cd$ is \emph{perfect} (respectively, \emph{superperfect}) if for every countable (respectively, small) family of morphisms $M_{i}\ra N_{i}\ko i\in I$, of $\cd$ such that the coproducts $\coprod_{I}M_{i}$ and $\coprod_{I}N_{i}$ exist, the induced map\[\cd(P,\coprod_{I}M_{i})\ra\cd(P,\coprod_{I}N_{i})
\]
is surjective provided every map
\[\cd(P,M_{i})\ra\cd(P,N_{i})\ko i\in I
\]
is surjective. Particular cases of superperfect objects are \emph{compact} objects, \ie objects $P$ such that the functor $\cd(P,?)$ preserves small coproducts.

\subsection{Milnor colimits}\label{Milnor colimit}

Now we recall a crucial construction which formally imitates the construction of the direct limit in an abelian category.
Let $\cd$ be a triangulated category and let
\[M_{0}\arr{f_{0}}M_{1}\arr{f_{1}}M_{2}\arr{f_{2}}\dots
\]
be a sequence of morphisms of $\cd$ such that the coproduct $\coprod_{n\geq 0}M_{n}$ exists in $\cd$. The \emph{Milnor colimit}\index{Milnor!colimit} of this sequence, denoted by $\Mcolim M_{n}$, is given, up to non-unique isomorphism, by the triangle\[\coprod_{n\geq 0}M_{n}\arr{\id-\sigma}\coprod_{n\geq 0}M_{n}\arr{\pi} \Mcolim M_{n}\ra\coprod_{n\geq 0}M_{n}[1],\]
where the morphism $\sigma$ has components
\[M_{n}\arr{f_{n}} M_{n+1}\arr{\can}\coprod_{p\geq 0}M_{p}.
\]
The above triangle is the \emph{Milnor triangle} (\cf \cite{Milnor1962, Keller1998b}) associated to the sequence $f_{n}\ko n\geq 0$. The notion of Milnor colimit has appeared in the literature under the name of \emph{homotopy colimit} (\cf \cite[Definition 2.1]{BokstedtNeeman1993}, \cite[Definition 1.6.4]{Neeman2001}) and \emph{homotopy limit} (\cf \cite[subsection 5.1]{Keller1994a}). However, we think it is better to keep this terminology for the notions appearing in the theory of derivators \cite{Maltsiniotis2001, Maltsiniotis2005, CisinskiNeeman2005} and in the theory of model categories \cite{Hirschhorn2003}.
\subsection{Generation of triangulated categories}\label{Generation of triangulated categories}
Let us consider three ways in which a triangulated category $\cd$ can be generated by a class $\cq$ of objects:
\begin{enumerate}[1)]
\item $\cd$ is \emph{generated} by $\cq$ if an object $M$ of $\cd$ is zero whenever\[\cd(Q[n],M)=0
\]
for every object $Q$ of $\cq$ and every integer $n\in\Z$. In this case, we say that $\cq$ is a \emph{class of generators} of $\cd$ and that $\cq$ \emph{generates} $\cd$. A triangulated category with small coproducts is \emph{compactly generated} if it is generated by a set of compact objects.\item $\cd$ satisfies the \emph{principle of infinite d\'{e}vissage} with respect to $\cq$ if $\cd=\Tria_{\cd}(\cq)$. In this situation, $\cq$ generates $\cd$.\item $\cd$ is \emph{exhaustively generated} by $\cq$ if the following conditions hold:\begin{enumerate}[3.1)]
\item Small coproducts of objects of $\bigcup_{m\geq 0}\Sum(\cq)^{*m}$ exist in $\cd$.\item For each object $M$ of $\cd$ there exists an integer $i\in\Z$ and a triangle\[\coprod_{n\geq 0}Q_{n}\ra\coprod_{n\geq 0}Q_{n}\ra M[i]\ra \coprod_{n\geq 0}Q_{n}[1]\]
in $\cd$ with $Q_{n}\in\bigcup_{m\geq 0}\Sum(\cq)^{*m}$.
\end{enumerate}
Notice that, in this situation, $\cd$ satisfies the principle of infinite dévissage with respect to $\cq$. If $\cq=\cp^+$ for some set $\cp$, then we also say that $\cd$ is \emph{exhaustively generated to the left} by $\cp$.\end{enumerate}

The following are two examples of exhaustively generated triangulated categories:
\begin{example}\label{Krause on perfect}
Let $\cd$ be a triangulated category with small coproducts, and let $\cp$ be a set of objects of $\cd$ which are perfect in $\Tria(\cp)$. As proved by \cite[Theorem A]{Krause2002}, every object of $\Tria(\cp)$ is the Milnor colimit of a sequence\[P_{0}\arr{f_{0}}P_{1}\arr{f_{1}}P_{2}\ra\dots
\]
of morphisms of $\cd$ where $P_{n}$ is an $n$th extension of small coproducts of shifts of objects of $\cp$. This shows that $\Tria(\cp)$ is exhaustively generated by the set formed by all the shifts of objects of $\cp$. In particular, the derived category $\cd\ca$ of a small dg category $\ca$ is exhaustively generated by all the shifts of the representable modules $A^{\we}:=\ca(?,A)\ko A\in\ca$.\end{example}

\begin{example}\label{Souto on perfect}
Let $\cd$ be a triangulated category with small coproducts, and let $\cp$ be a set of perfect objects of $\cd$. As proved in \cite[Theorem 2.2]{Souto2004}, we have that $\Susp(\cp)$ is an aisle in $\cd$ and every object of $\Susp(\cp)$ is a Milnor colimit of a sequence\[P_{0}\arr{f_{0}}P_{1}\arr{f_{1}}P_{2}\ra\dots
\]
of morphisms of $\cd$ where $P_{n}$ is an $n$th extension of small coproduct of non-negative shifts of objects of $\cp$. In particular, $\bigcup_{n\in\Z}\Susp(\cp)[n]$ is exhaustively generated to the left by $\cp$.\end{example}

\subsection{The right bounded derived category of a dg category}

Let $\ca$ be a small dg category. Since the representable dg right $\ca$-modules $A^{\we}\ko A\in\ca$, are compact objects of the derived category $\cd\ca$ of $\ca$, then $\Susp(\{A^{\we}\}_{A\in\ca})$ is an aisle in $\cd\ca$, which will be denoted by $\cd^{\leq 0}\ca$. Its associated coaisle, denoted by $\cd^{>0}\ca$, consists of those modules $M$ with cohomology concentrated in positive degrees, \ie $\H nM(A)=0$ for each $A\in\ca$ and $n\leq 0$. For each integer $n\in\Z$ we put\[\cd^{\leq n}\ca:=\cd^{\leq 0}\ca[-n]
\]
and
\[\cd^{>n}\ca:=\cd^{>0}[-n],
\]
and denote by $\tau^{\leq n}$ and $\tau^{>n}$ the torsion and torsionfree functors, respectively, corresponding to the t-structure $(\cd^{\leq n}\ca,\cd^{>n}\ca)$.
The following lemma ensures that, in case the dg category $\ca$ has cohomology concentrated in non-positive degrees, the aisle $\cd^{\leq n}\ca$ admits a familiar description in terms of cohomology.
\begin{lemma}\label{looking for coherence}
Let $\ca$ be a small dg category with cohomology concentrated in degrees $(-\infty,m]$ for some integer $m\in\Z$. For a dg $\ca$-module $M$ we consider the following assertions:\begin{enumerate}[1)]
\item $M\in\cd^{\leq s}\ca$.
\item $\H iM(A)=0$ for each integer $i>m+s$ an every object $A$ of $\ca$.\end{enumerate}
Then $1)\Rightarrow 2)$ and, in case $m=0$, we also have $2)\Rightarrow 1)$.\end{lemma}
\begin{proof}
$1)\Rightarrow 2)$ Since $M[s]$ belongs to $\Susp(\{A^{\we}\}_{A\in\ca})$, there exists a triangle in $\cd\ca$\[\coprod_{n\geq 0}P_{n}\ra\coprod_{n\geq 0}P_{n}\ra M[s]\ra \coprod_{n\geq 0}P_{n}[1]\]
with $P_{n}\in\Sum(\{A^{\we}\}^{+}_{A\in\ca})^{*n}$ for each $n\geq 0$ (\cf for instance Example \ref{Souto on perfect}). Then, for each $A\in\ca$ we get the long exact sequence of cohomology\[\dots\ra\coprod_{n\geq 0}\H iP_{n}(A)\ra\H {i+s}M(A)\ra\coprod_{n\geq 0}\H {i+1}P_{n}(A)\ra\dots\]
with $\H iP_{n}(A)\cong(\cd\ca)(A^{\we},P_{n}[i])=0$ for each $i>m$.

$2)\Rightarrow 1)$ Consider the triangle in $\cd\ca$
\[M'\ra M\ra M''\ra M'[1]
\]
with $M'\in\cd^{\leq s}\ca$ and $M''\in(\cd^{\leq s}\ca)^{\bot}$. In particular, $\H {i}M''(A)=0$ for each $A\in\ca$ and each $i\leq s$. The aim is to prove that $\H iM''(A)=0$ for each $A\in\ca$ and each $i\in\Z$. Thus, consider the long exact sequence of cohomology\[\dots\ra\H iM(A)\ra\H iM''(A)\ra\H {i+1}M'(A)\ra\dots
\]
By using 1) and the extra assumption on $\ca$, we have that $\H {i}M'(A)=0$ for each $i>s$ and, by hypothesis, $\H iM(A)=0$ for each $i>s$. This implies that $\H iM''(A)=0$ for each $i>s$.\end{proof}

For an arbitrary small dg category $\ca$, the t-structure $(\cd^{\leq 0}\ca,\cd^{>0}\ca)$ is said to be the \emph{canonical t-structure} on $\cd\ca$. We will write\[\cd^-\ca:=\bigcup_{n\in\Z}\cd^{\leq n}\ca,
\]
and we will refer to $\cd^-\ca$ as the \emph{right bounded derived category} of $\ca$. These names are justified by the Lemma \ref{looking for coherence}.
\begin{remark}
Notice that $\cd^-\ca$ is not closed under small coproducts in $\cd\ca$. Indeed, given $A\in\ca$, the coproduct $\coprod_{n\in\Z}A^{\we}[n]$ does not belong to $\cd^-\ca$. Also, notice that $\cd^-\ca$ is exhaustively generated to the left by the free $\ca$-modules $A^{\we}\ko A\in\ca$.\end{remark}


\section{Lifting of TTF triples}\label{Lifting of TTF triples}

\subsection{General criterion}\label{General criterion}

\begin{definition}\label{definition of restriction}
Let $\cd$ be a triangulated category and let $\cd'$ be a full triangulated subcategory of $\cd$. We say that a TTF triple $(\cx,\cy,\cz)$ on $\cd$\emph{restricts to} or \emph{is a lifting of} a TTF triple $(\cx',\cy',\cz')$ on $\cd'$ if we have\[(\cx\cap\cd',\cy\cap\cd',\cz\cap\cd')=(\cx',\cy',\cz').
\]
That is to say, $\cx'$ is the full subcategory of $\cd'$ formed by those objects of $\cd'$ which are in $\cx$, and analogously with the other subcategories.In this case, we say that $(\cx',\cy',\cz')$ \emph{lifts to} or \emph{is the restriction of}  $(\cx,\cy,\cz)$.\end{definition}

\begin{definition}
A class $\cp$ of objects of a triangulated category $\cd$ is \emph{recollement-defining} if the class $\cy$ of those objects which are right orthogonal to all the shifts  of objects of $\cp$ is both an aisle and a coaisle in $\cd$, \ie $\cy$ fits in a TTF triple $(^{\bot}\cy,\cy,\cy^{\bot})$ on $\cd$.\end{definition}

\begin{proposition}\label{restriction}
Let $\cd$ be a triangulated category with small coproducts and let
$\cd'$ be a full triangulated subcategory containing a set $\cq$
of generators of $\cd$. For a TTF triple $(\cx',\cy',\cz')$ on
$\cd'$ the following assertions are equivalent:
\begin{enumerate}[1)]
\item $(\cx',\cy',\cz')$ is the restriction of a TTF triple on $\cd$.
\item There is a set $\cp$ of objects of $\cx'$ such that:
\begin{enumerate}[2.1)]
\item $\cp$ is recollement-defining in $\cd$.
\item If an object of $\cd$ is right orthogonal to all the shifts of objects of $\cp$, then it is right orthogonal to all the objects of $\cx'$.\end{enumerate}
\item The objects of $\cx'$ form a recollement-defining class of $\cd$.
\end{enumerate}
Moreover, we can take $\cp=(\tau_{\cx'}z'\tau^{\cz'})(\cq)$.
\end{proposition}
\begin{proof}
$1)\Rightarrow 2)$ Let $(\cx,\cy,\cz)$ be a TTF triple on $\cd$ which restricts to $(\cx',\cy',\cz')$, and let $\cq$ be a set of generators of $\cd$ contained in $\cd'$. Notice that, for each object $Q$ of $\cq$, the torsion triangle associated to the t-structure $(\cy,\cz)$ can be taken to be\[\tau_{\cy'}(Q)\ra Q\ra \tau^{\cz'}(Q)\ra \tau_{\cy'}(Q)[1].
\]
Then, it is straightforward to check that $\tau^{\cz'}(\cq)$ is a set of generators of $\cz$. Since the composition\[\cz\arr{z}\cd\arr{\tau_{\cx}}\cx
\]
is a triangle equivalence, we have that $\cp:=(\tau_{\cx}z\tau^{\cz'})(\cq)$ is a set of generators of $\cx$. But, since $\tau^{\cz'}(\cq)$ is contained in $\cd'$, then we have $\cp=(\tau_{\cx'}z'\tau^{\cz'})(\cq)$, which is contained in $\cx'$. The fact that $(\cx,\cy)$ is a t-structure on $\cd$ implies that $\cy$ is the set of objects of $\cd$ which are right orthogonal to all the shifts of objects of $\cp$, and so $\cp$ is recollement-defining in $\cd$. Finally, the inclusions $\cy\subseteq\cx^{\bot}\subseteq\cx'^{\bot}$ prove 2.2).
$2)\Rightarrow 3)$ is clear.

$3)\Rightarrow 1)$ Consider $(\cx,\cy,\cz):=(^{\bot}(\cx'^{\bot}),\cx'^{\bot},(\cx'^{\bot})^{\bot})$, with orthogonals taken in $\cd$, which is a TTF triple on $\cd$. Since $(\cx',\cy',\cz')$ is a TTF triple on $\cd'$, then we have $\cy'=\cx'^{\bot}\cap\cd'=\cy\cap\cd'$. Let us prove now $\cx'=\cx\cap\cd'$. The inclusion $\subseteq$ is clear. Conversely, let $X$ be an object of $\cx\cap\cd'$ and consider the triangle\[\tau_{\cx'}(X)\ra X\ra \tau^{\cy'}(X)\ra \tau_{\cx'}(X)[1].
\]
Its two terms on the left belong to $\cx$. Then $\tau^{\cy'}(X)\in \cx\cap\cy'\subseteq\cx\cap\cy=\{0\}$ and so $X\in\cx'$. Now, we have the following inclusions\[\cz\cap\cd'=\cy^{\bot}\cap\cd'\subseteq\cy'^{\bot}\cap\cd'=\cz'.
\]
Finally, let $\cq$ be the set of generators of $\cd$ contained in $\cd'$. It is easy to prove that $\tau^{\cy}(\cq)^{\bot}=\cz$. Also, notice that $\tau^{\cy}(\cq)\subseteq \cy\cap\cd'=\cy'$. Therefore,\[\cz'=\cy'^{\bot}\cap\cd'\subseteq \tau^{\cy}(\cq)^{\bot}\cap\cd'=\cz\cap\cd'.\]
\end{proof}

\begin{corollary}\label{bijection between restricts and are restricted}
Under the hypotheses of Proposition \ref{restriction}, the map
\[(\cx,\cy,\cz)\mapsto(\cx\cap\cd',\cy\cap\cd',\cz\cap\cd')
\]
defines a bijection between:
\begin{enumerate}[1)]
\item TTF triples on $\cd$ which restricts to TTF triples on $\cd'$.
\item TTF triples on $\cd'$ which are restriction of TTF triples on $\cd$.\end{enumerate}
\end{corollary}
\begin{proof}
Of course, the map is surjective. Now, let $\cq\subseteq\cd'$ a set of generators of $\cd$ and let $(\cx,\cy,\cz)$ be a TTF triple such that $(\cx',\cy',\cz')=(\cx\cap\cd',\cy\cap\cd',\cz\cap\cd')$ is a TTF triple on $\cd'$. Then, the proof of Proposition \ref{restriction} shows that $\cy$ is precisely the class of objects of $\cd$ which are right orthogonal to all the shifts of objects of $(\tau_{\cx'}z'\tau^{\cz'})(\cq)$. This implies the injectivity.\end{proof}

\subsection{`Right bounded' triangulated subcategories}\label{`Right bounded' triangulated subcategories}
Let $\cq$ be a set of objects of a triangulated category $\cd$ with small coproducts. Let us assume that $\Susp(\cq)$ is an aisle in $\cd$. This is the case, for instance if the objects of $\cq$ are perfect (\cf \cite{Krause2002, Souto2004}). Notice that, in case $\Susp(\cq)$ is an aisle in $\cd$, then $\Susp(\cq)=\aisle(\cq)$, \ie $\Susp(\cq)$ is the smallest aisle in $\cd$ containing $\cq$.
We are interested in the interplay between TTF triples on abstract `unbounded' triangulated categories and TTF triples on abstract `right bounded triangulated' categories. More precisely, we are interested in the interplay between TTF triples on $\cd$ and TTF triples on the full triangulated subcategory $\cd':=\bigcup_{n\in\Z}\aisle(\cq)[n]$ of $\cd$. A good example to keep in mind is $\cd=\cd\ca$ and $\cd'=\cd^-\ca$ for a small dg category $\ca$.
First we need to understand better the interplay between $\cd$ and $\cd'$.
\begin{lemma}\label{preservation of coproducts and compacity}
The following assertions hold:
\begin{enumerate}[1)]
\item The inclusion functor $\iota:\cd'\hookrightarrow\cd$ preserves small coproducts.\item If $\Susp(\cq)^{\bot}$ is closed under small coproducts, then an object $P$ of $\cd'$ is compact (respectively, perfect, superperfect) in $\cd'$ if and only if it is compact (respectively, perfect, superperfect) in $\cd$.\end{enumerate}
\end{lemma}
\begin{proof}
1) Let $D'_{i}\ko i\in I$, be a family of objects of $\cd'$ whose coproduct exists in $\cd'$. We write $\coprod_{i\in I}D'_{i}$ for the coproduct in $\cd$, $D'$ for the coproduct in $\cd'$ and $v_{i}:D'_{i}\ra D'$ for the canonical morphisms. For simplicity, put $\aisle(\cq)[k]=\cu_{k}$. Therefore, we have a chain\[\dots\subseteq\cu_{k+1}\subseteq\cu_{k}\subseteq\cu_{k-1}\subseteq\dots\subseteq\cd'\]
of aisles in $\cd$ whose union is $\cd'$.

\emph{Claim: If $m\ko n\in\Z$ are integers such that $D'\in\cu_{n}$ and $D'_{i}\in\cu_{m}\setminus\cu_{m+1}$ for some $i\in I$, then $n\leq m$.} Indeed, fix such an $i$ and assume $n>m$ and consider the triangle\[\tau_{\cu_{n}}(D'_{i})\ra D'_{i}\arr{f}\tau^{\cu_{n}^{\bot}}(D'_{i})\ra \tau_{\cu_{n}}(D'_{i})[1].\]
Since the two first vertices of this triangle belong to $\cd'$, then so does $\tau^{\cu_{n}^{\bot}}(D'_{i})$. Hence, by using the universal property of the coproduct, we have that $f$ induces a morphism\[\tilde{f}:D'\ra \tau^{\cu_{n}^{\bot}}(D'_{i})
\]
such that
\[\tilde{f}v_{j}=\begin{cases}f & \text{ if }j=i, \\
0 & \text{ otherwise.}\end{cases}
\]
Since $D'\in\cu_{n}$, then $\tilde{f}=0$ and so $f=0$. Therefore, $D'_{i}$ is a direct summand of $\tau_{\cu_{n}}(D'_{i})$. This implies that $D'_{i}$ belongs to $\cu_{n}$, and so it belongs to $\cu_{m+1}$, which is a contradiction.
Consider the following two situations:

\emph{First situation:} For each $i\in I$ we have $D'_{i}\in\bigcap_{k\in\Z}\cu_{k}$. Since aisles are closed under small coproducts, this implies that the coproduct $\coprod_{i\in I}D'_{i}$ belongs to $\bigcap_{k\in\Z}\cu_{k}$, and so to $\cd'$. Hence $D'\cong\coprod_{i\in I}D'_{i}$.
\emph{Second situation:} There exists $j\in I$ such that $D'_{j}\in\cu_{m}\setminus\cu_{m+1}$. Given $i\in I$, put $m_{i}$ for the maximum of the set of those integers $k\in\Z$ such that $D'_{i}\in\cu_{k}$. Put $m_{i}=\infty$ if $D'_{i}\in\bigcap_{k\in\Z}\cu_{k}$. Thanks to the claim, we know that, in any case, $m_{i}\geq n$ for each $i\in I$. Then $D'_{i}\in\cu_{n}$ for every $i\in I$, and so $\coprod_{i\in I}D'_{i}\in\cu_{n}$. Again, this implies $\coprod_{i\in I}D'_{i}\cong D'$.
2) Assertion 1) implies that if $P\in\cd'$ is compact in $\cd$ then it is also compact in $\cd'$. Conversely, let $P\in\cd'$ be compact in $\cd'$ and fix an integer $n\in\Z$ such that $P\in\cu_{n}$. If $D_{i}\ko i\in I$, is a family of objects of $\cd$, then we have isomorphisms\[\cd(P,D_{i})\cong\cu_{n}(P,\tau_{\cu_{n}}(D_{i}))=\cd'(P,\tau_{\cu_{n}}(D_{i}))\]
for each $i\in I$, and
\[\cd(P,\coprod_{i\in I}D_{i})\cong\cu_{n}(P,\tau_{\cu_{n}}(\coprod_{i\in I}D_{i}))=\cd'(P,\tau_{\cu_{n}}(\coprod_{i\in I}D_{i})).\]
By hypothesis, $\cu^{\bot}_{n}$ is closed under small coproducts. This is equivalent to the fact that $\tau_{\cu_{n}}$ preserves small coproducts, and so we have a canonical isomorphism\[\coprod_{i\in I}\tau_{\cu_{n}}(D_{i})\arr{\sim}\tau_{\cu_{n}}(\coprod_{i\in I}D_{i}).\]
Finally, we have the commutative diagram
\[\xymatrix{\coprod_{i\in I}\cd(P,D_{i})\ar[r]_{\sim\ \ \ \ \ \ }\ar[dd]^{\can} & \coprod_{i\in I}\cd'(P,\tau_{\cu_{n}}(D_{i}))\ar[d]_{\wr}^{\can} \\&\cd'(P,\coprod_{i\in I}\tau_{\cu_{n}}(D_{i}))\ar[d]_{\wr}^{\can} \\
\cd(P,\coprod_{i\in I}D_{i})\ar[r]_{\sim\ \ \ \ \ \ } & \cd'(P,\tau_{\cu_{n}}(\coprod_{i\in I}D_{i}))}
\]
where the morphisms `can' are the canonical ones. This proves that $P$ is compact in $\cd$. The case of $P$ being (super)perfect follows similarly using adjunction.\end{proof}

\begin{proposition}\label{easy restriction}
Assume that $\cq$ is a set of perfect generators of $\cd$ such
that  $\aisle(\cq)^{\bot}$ is closed under small coproducts. Let
$(\cx',\cy',\cz')$ be a TTF triple on $\cd'$ such that $\cx'$ is
exhaustively generated to the left by a set $\cp$ whose objects
are superperfect in $\cx'$. Then,
\begin{enumerate}[1)]
\item The objects of $\cp$ are superperfect in $\cd'$.
\item $(\Tria_{\cd}(\cp), \Tria_{\cd}(\cp)^{\bot}, (\Tria_{\cd}(\cp)^{\bot})^{\bot})$ is a TTF triple on $\cd$ which restricts to $(\cx',\cy',\cz')$.\end{enumerate}
\end{proposition}
\begin{proof}
1) Let $P$ be an object of $\cp$. Let $\alpha_{i}:M_{i}\ra N_{i}\ko i\in I$, be a family of morphism of $\cd'$ such that the induced maps $\cd'(P,M_{i})\ra\cd'(P,N_{i})\ko i\in I$, are surjective. In other words, the maps $\cx'(P,\tau_{\cx'}M_{i})\ra\cx'(P,\tau_{\cx'}N_{i})\ko i\in I$, are surjective. Assume the coproducts $\coprod_{I}M_{i}$ and $\coprod_{I}N_{i}$ exist in $\cd'$. We have to prove that the induced map $\cd'(P,\coprod_{I}M_{i})\ra\cd'(P,\coprod_{I}N_{i})$ is surjective. For each $i\in I$ we consider the triangle\[x'\tau_{\cx'}M_{i}\arr{f_{i}}M_{i}\arr{g_{i}}y'\tau^{\cy'}M_{i}\arr{h_{i}}x'\tau_{\cx'}M_{i}[1]\]
of $\cd'$ associated to the t-structure $(\cx',\cy')$. Since both $\tau^{\cy'}$ and $y'$ preserve small coproducts, the coproduct $\coprod_{I}y'\tau^{\cy'}M_{i}$ exists in $\cd'$ and the canonical morphism $\coprod_{I}y'\tau^{\cy'}M_{i}\ra y'\tau^{\cy'}\coprod_{I}M_{i}$ is an isomorphism. The existence of the coproducts $\coprod_{I}M_{i}$ and $\coprod_{I}y'\tau^{\cy'}M_{i}$ implies that the coproduct $\coprod_{I}x'\tau_{\cx'}M_{i}$ exists in $\cd'$ and that\[\coprod_{I}x'\tau_{\cx'}M_{i}\arr{\coprod_{I}f_{i}}\coprod_{I}M_{i}\arr{\coprod_{I}g_{i}}\coprod_{I}y'\tau^{\cy'}M_{i}\arr{\coprod_{I}h_{i}}x'\tau_{\cx'}\coprod_{I}M_{i}[1]\]
is a triangle of $\cd'$. Hence, the canonical morphism $\coprod_{I}x'\tau_{\cx'}M_{i}\ra x'\tau_{\cx'}\coprod_{I}M_{i}$ is an isomorphism. Of course, we can proceed similarly with the objects $N_{i}\ko i\in I$. Then, the map $\cd'(P,\coprod_{I}M_{i})\ra\cd'(P,\coprod_{I}N_{i})$ is isomorphic to $\cx'(P,\coprod_{I}\tau_{\cx'}M_{i})\ra\cx'(P,\coprod_{I}\tau_{\cx'}N_{i})$, which is surjective since $P$ is superperfect in $\cx'$.
2) Lemma \ref{preservation of coproducts and compacity} implies that the objects of $\cp$ are also superperfect in $\cd$. Then, by using Brown representability theorem \cite{Krause2002} we deduce that $\cx:=\Tria_{\cd}(\cp)$ is an aisle in $\cd$. Notice that the corresponding coaisle $\cy:=\Tria_{\cd}(\cp)^{\bot}$ is formed by those objects which are right orthogonal to all the shifts of objects of $\cp$, and so it is closed under small coproducts. Since $\tau^{^\cy}(\cq)$ is a set of perfect generators of $\cy$, we can use again Brown representability theorem to deduce that $\cy$ is an aisle in $\cd$. Put $\cz:=\cy^{\bot}$. Of course, condition 2.1) of Proposition \ref{restriction} is satisfied. Let us prove that condition 2.2) of this proposition also holds. For this, first notice that thanks to Lemma \ref{preservation of coproducts and compacity}, we know that the inclusion functor $\cx'\hookrightarrow\cd$ preserves small coproducts for it is the composition of the coproduct-preserving inclusions $\cx'\hookrightarrow\cd'\hookrightarrow\cd$. Now, for every object $X\in\cx'$ there exists an integer $i\in\Z$ such that $X$ fits into a triangle of $\cx'$ (and so of $\cd$)\[\coprod_{n\geq 0}P_{n}[i]\ra \coprod_{n\geq 0}P_{n}[i]\ra X\ra \coprod_{n\geq 0}P_{n}[i+1]\]
with $P_{n}\in\bigcup_{m\geq 0}\Sum(\cp^+)^{*m}$ for each $n\geq 0$. Let $M$ be an object of $\cd$ which is right orthogonal to all the shifts of objects of $\cp$. Then, we get a long exact sequence\[\dots\ra\cd(\coprod_{n\geq 0}P_{n}[i+1],M)\ra\cd(X,M)\ra\cd(\coprod_{n\geq 0}P_{n}[i],M)\ra\dots\]
in which
\[\cd(\coprod_{n\geq 0}P_{n}[i+1],M)=\cd(\coprod_{n\geq 0}P_{n}[i],M)=0,\]
and so $\cd(X,M)=0$.

\end{proof}

\section{Restriction of TTF triples}\label{Restriction of TTF triples}

\subsection{General criterion}\label{General criterion bis}

The general criterion to restrict t-structures is the following well-known lemma, which already appeared in the work of A. A. Beilinson, J. Bernstein and P. Deligne (\cf \cite[paragraph 1.3.19]{BeilinsonBernsteinDeligne}):
\begin{lemma}\label{general restriction of t-structures}
Let $(\cu,\cv[1])$ be a t-structure on a triangulated category $\cd$, and let $\cd'$ be a strictly(=closed under isomorphisms) full triangulated subcategory of $\cd$. The following assertions are equivalent:\begin{enumerate}[1)]
\item $(\cd'\cap\cu,\cd'\cap\cv[1])$ is a t-structure on $\cd'$.
\item $(u\tau_{\cu})(\cd')\subseteq\cd'$.
\end{enumerate}
\end{lemma}

\subsection{`Right bounded' triangulated subcategories}\label{`Right bounded' triangulated subcategories bis}
We present now a very particular situation in which condition 2) of the lemma above can be improved.
As in subsection \ref{`Right bounded' triangulated subcategories}, let $\cd$ be a triangulated category with small coproducts, and let $\cq$ be a set of objects of  $\cd$ such that $\Susp(\cq)$ is an aisle in $\cd$ (and so $\Susp(\cq)=\aisle(\cq)$). Let $\cd':=\bigcup_{n\in\Z}\aisle(\cq)[n]$.
\begin{proposition}\label{particular restriction of t-structures}
Assume that $\cq$ is a set of perfect generators of $\cd$ such
that $\aisle(\cq)^{\bot}$ is closed under small coproducts. Let
$(\cu,\cv)$ be a t-structure on $\cd$ such that $\cu$ is
triangulated and $\cv$ is closed under coproducts. The following
assertions are equivalent:
\begin{enumerate}[1)]
\item $(\cd'\cap\cu,\cd'\cap\cv)$ is a t-structure on $\cd'$.
\item $(u\tau_{\cu})(\cq)\subseteq\aisle(\cq)[n]$ for some integer $n$.
\end{enumerate}
\end{proposition}
\begin{proof}
$1)\Rightarrow 2)$ Thanks to Lemma \ref{general restriction of t-structures}, it suffices to prove that $(u\tau_{\cu})(\cd')\subseteq\cd'$ implies condition 2) of the proposition. Since $N:=\coprod_{Q\in\cq}Q$ belongs to $\aisle(\cq)$, there exists an integer $n$ such that $u\tau_{\cu}N$ belongs to $\aisle(\cq)[n]$. Now notice that for each $Q\in\cq$ we have that $u\tau_{\cu}Q$ is a direct summand of $u\tau_{\cu}N$. But since $\aisle(\cq)[n]$ is closed under Milnor colimits in $\cd$ then it is also closed under direct summands. This implies that $u\tau_{\cu}Q$ belongs to $\aisle(\cq)[n]$.
$2)\Rightarrow 1)$ Thanks to Lemma \ref{general restriction of t-structures}, it suffices to prove the inclusion $(u\tau_{\cu})(\cd')\subseteq\cd'$. Let $N$ be an object of $\cd'$ and fix an integer $i$ such that $N[i]$ belongs to $\aisle(Q)$. The proof of \cite[Theorem 2.2]{Souto2004} shows us that $N[i]$ is the Milnor colimit of a sequence\[M_{0}\arr{f_{0}}M_{1}\arr{f_{1}}M_{2}\arr{f_{2}}\dots
\]
where $M_{n}\in\Sum(\cq^+)^{*n}$ for each $n\geq 0$. Now, since $u\tau_{\cu}$ commutes with small coproducts, by applying it to the corresponding Milnor triangle we get that $u\tau_{\cu}N[i]$ belongs to $\cd'$, and then so does $u\tau_{\cu}N$.\end{proof}

\begin{remark}
If in Proposition \ref{particular restriction of t-structures} the set $\cq$ is finite, then one can replace condition 2) by: $(u\tau_{\cu})(\cq)\subseteq\cd'$.\end{remark}

\begin{corollary}\label{conditions to restrict}
Assume that $\cq$ is a set of perfect generators of $\cd$ such
that $\aisle(\cq)^{\bot}$ is closed under small coproducts. The
following assertions are equivalent for a TTF triple
$(\cx,\cy,\cz)$ on $\cd$:
\begin{enumerate}[1)]
\item $(\cd'\cap\cx,\cd'\cap\cy,\cd'\cap\cz)$ is a TTF triple on $\cd'$.
\item The following conditions hold:
\begin{enumerate}[2.1)]
\item $(x\tau_{\cx})(\cq)\subseteq\aisle(\cq)[n]$ for some integer $n$.
\item $(y\tau_{\cy})(\cd')\subseteq\cd'$.
\end{enumerate}
\end{enumerate}
\end{corollary}

\begin{example}\label{restriction in examples}
Let $I$ be a two-sided ideal of a $k$-algebra $A$, and assume the canonical projection $\pi:A\ra A/I$ is a \emph{homological epimorphism} in the sense of W. Geigle and H. Lenzing \cite{GeigleLenzing}. We know (\cf \cite[Example 5.3.4]{NicolasTesis}) that in this case $\cd A$ is a recollement of $\cd(A/I)$ and $\Tria_{\cd A}(I)$:\[\xymatrix{\cd(A/I)\ar[rr]^{\pi_{*}} && \cd A\ar@/_1.2pc/[ll]_{?\otimes^{\L}_{A}A/I}\ar@/_-1pc/[ll]^{\R Hom_{A}(A/I,?)}\ar[rr]^{?\otimes^{\L}_{A}I} && \Tria_{\cd A}(I)\ar@/_1.2pc/[ll]\ar@/_-1pc/[ll]^x}
\]
where $x$ is the inclusion functor. Let $\cc_{dg}A$ be the dg category whose objects are the complexes of $A$-modules and whose morphisms are given by complexes of $k$-modules $\cc_{dg}(A)(L,M)$, with $n$th component formed by the morphisms of $\Z$-graded $k$-modules homogeneous of degree $n$ and with differential given by the commutator $d(f)=d_{M}f-(-1)^{|f|}fd_{L}$, where $|f|$ is the degree of $f$. Notice that the corresponding category of $0$-cocycles $\Zy 0(\cc_{dg}A)$ is the category $\cc A$ of complexes of $A$-modules and the corresponding category of $0$-cohomology $\H 0(\cc_{dg}A)$ is the category $\ch A$ of complexes of $A$-modules  up to homotopy. In case $I$ is compact in $\Tria_{\cd A}(I)$, the proof of \cite[Theorem 4.3]{Keller1994a} implies that $\cd A$ is a recollement of $\cd(A/I)$ and $\cd C$, where $C$ is the dg algebra $(\cc_{dg}A)(\textbf{i}I,\textbf{i}I)$ and $\textbf{i}:\cd A\ra\ch A$ is the \emph{fibrant replacement functor} (\cf \cite{Keller2006b}). Indeed, the dg $A$-$C$-bimodule $\textbf{i}I$ induces mutually quasi-inverse triangle equivalences\[\xymatrix{\Tria_{\cd A}(I)\ar@<1ex>[rr]^{\R\Hom_{A}(\textbf{i}I,?)} && \cd C.\ar@<1ex>[ll]^{?\otimes^{\L}_{C}\textbf{i}I}}
\]
Thanks to Corollary \ref{conditions to restrict}, we know that the associated TTF triple restricts to $\cd^-A$ if and only if the following conditions hold:\begin{enumerate}[1)]
\item $A\otimes^{\bf L}_{A}I\cong I$ belongs to $\cd^-A$,
\item $\textbf{R}\Hom_{A}(A/I,M)$ belongs to $\cd^-(A/I)$ for each $M$ in $\cd^-A$.\end{enumerate}
Of course, the first condition always holds. Thanks to S.~K\"{o}nig's criterion explained at the begining of the proof of \cite[Theorem 1]{Konig1991}, we have that the second condition holds if and only if $A/I$ has finite projective dimension regarded as a right $A$-module or, equivalently, $I$ has finite projective dimension regarded as a right $A$-module. Assume then that $I_{A}$ has finite projective dimension and also that it is compact in $\Tria_{\cd A}(I)$. In this case the mutually quasi-inverse triangle equivalences between $\Tria_{\cd A}(I)$ and $\cd C$restrict to mutually quasi-inverse triangle equivalences
\[\xymatrix{\Tria_{\cd A}(I)\cap\cd^-A\ar@<1ex>[r] & \cd^-C.\ar@<1ex>[l]
}
\]
Therefore, $\cd^-A$ is a recollement of $\cd^-(A/I)$ and $\cd^-C$. This example contains as particular cases the recollement data of Corollary 11, Corollary 12 and Corollary 15 of \cite{Konig1991}, and describes functors appearing in those recollement data as restrictions of total derived functors.\end{example}

\section{Recollement of right bounded derived categories}\label{Recollements of right bounded derived categories}
All through this section the appearing dg categories are small.

\subsection{Bounds}\label{bounds}

\begin{definition}\label{fibrant replacement for sets}
Let $\ca$ be a dg category. Consider the corresponding dg category $\cc_{dg}\ca$ (\cf \cite{Keller2006b}), which is the `dg generalization to several objects' of the dg category $\cc_{dg}A$ associated to an algebra $A$ appearing in Example \ref{restriction in examples}. A \emph{fibrant replacement} of a set $\cp$ of objects of the derived category $\cd\ca$ is a full subcategory $\cb$ of $\cc_{dg}\ca$ formed by the fibrant replacements $\textbf{i}P$, in the sense of \cite{Keller2006b}, of the modules $P$ of $\cp$.\end{definition}

Notice that $\cb$ is a dg category and we have a dg $\cb$-$\ca$-bimodule $X$ defined by $X(A,B):=B(A)$ for $A$ in $\ca$ and $B$ in $\cb$. It is well-known (\cf \cite{Keller1994a,Keller2006b}) that this gives rise to a funtor\[\ch om_{\ca}(X,?):\cc_{dg}\ca\ra\cc_{dg}\cb
\]
which induces triangle functors
\[\ch om_{\ca}(X,?):\ch\ca\ra\ch\cb
\]
and
\[\R\ch om_{\ca}(X,?):\cd\ca\ra\cd\cb.
\]

\begin{definition}\label{dually right bounded}
Under the conditions above, we say that:
\begin{enumerate}[1)]
\item $\cp$ is \emph{right bounded} if $\cp\subseteq\cd^{\leq n}\ca$ for some $n\in\Z$.\item $\cp$ is \emph{dually right bounded} if the functor
\[\R\ch om_{\ca}(X,?):\cd\ca\ra\cd\cb
\]
sends an object of $\cd^-\ca$ to an object of $\cd^-\cb$.
\end{enumerate}
\end{definition}

\emph{A priori}, the notion of ``dually right bounded'' depends on
the fibrant replacement of $\cp$, however this is not really a
problem for our purposes. In the subsequent propositions we will
present the two situations in which we are most interested, where
the notion of ``dually right bounded'' is independent of the
fibrant replacement.

\begin{proposition}\label{characterization of dually right bounded}
Let $\ca$ be a dg category and $\cp$ a set of objects of $\cd^-\ca$. Assume that there exists an integer $m$ such that for every two objects $P$ and $P'$ of $\cp$ we have $(\cd\ca)(P,P'[i])=0$ for $i>m$. Consider the following assertions:\begin{enumerate}[1)]
\item $\cp$ is dually right bounded.
\item For each object $M$ of $\cd^-\ca$ there exists an integer $s_{M}$ such that $(\cd\ca)(P,M[i])=0$ for every $P\in\cp$ and every $i>m+s_{M}$.\end{enumerate}
Then 1) implies 2) and, if $m=0$, we also have that 2) implies 1).
\end{proposition}
\begin{proof}
Let $\cb$ be a fibrant replacement of the set $\cp$. Notice that the assumption on the set $\cp$ is equivalent to say that $\cb$ has cohomology concentrated in degrees $(-\infty,m]$. Let $X$ be the associated $\ca$-$\cb$-bimodule. Assertion 1) says that for each $M\in\cd^-\ca$ there exists an integer $s_{M}$ such that\[(\cc_{dg}\ca)(?,\textbf{i}M)_{|_{\cb}}\in\cd^{\leq s_{M}}\cb.
\]
Now, thanks to Lemma \ref{looking for coherence}, this implies that $(\cd\ca)(P,M[i])=0$ for every $P\in\cp$ and every $i>m+s_{M}$. As stated in Lemma \ref{looking for coherence}, in case $m=0$ we can go backward in the proof.\end{proof}

By using S.~K\"{o}nig's criterion which characterizes the bounded complexes of projective modules inside the right bounded derived category of an algebra (see the begining of the proof of \cite[Theorem 1]{Konig1991}), we deduce the following:
\begin{corollary}\label{dually right bounded are bounded complexes of projectives}Let $A$ be an ordinary algebra, and let $P$ an object of the right bounded derived category $\cd^-A$ of $A$ such that $(\cd A)(P,P[i])=0$ for $i\geq 1$. Then $P$ is dually right bounded if and only if it is quasi-isomorphic to a bounded complex of projective $A$-modules.\end{corollary}

\begin{proposition}\label{Morita theory for right bounded}
Let $\ca$ be a dg category and let $\cp$ be a set of objects of $\cd^-\ca$ such that:\begin{enumerate}[a)]
\item it is right bounded,
\item its objects are compact in $\Tria_{\cd\ca}(\cp)\cap\cd^-\ca$,
\item $\Tria_{\cd\ca}(\cp)\cap\cd^-\ca$ is exhaustively generated to the left by $\cp$.\end{enumerate}
Let $\cb$ be a fibrant replacement of $\cp$ and let $X$ be the associated $\cb$-$\ca$-bimodule. Then, the functor $?\otimes^{\L}_{\cb}X:\cd\cb\ra\cd\ca$ induces a triangle equivalence\[?\otimes^{\L}_{\cb}X:\cd^-\cb\arr{\sim}\Tria_{\cd\ca}(\cp)\cap\cd^-\ca,\]
and the following assertions are equivalent:
\begin{enumerate}[1)]
\item $\Tria_{\cd\ca}(\cp)\cap\cd^-\ca$ is an aisle in $\cd^-\ca$.
\item $\cp$ is dually right bounded.
\end{enumerate}
\end{proposition}
\begin{proof}
\emph{First step: The triangle functor
\[?\otimes^{\L}_{\cb}X:\cd\cb\ra\cd\ca
\]
induces a triangle functor
\[?\otimes^{\L}_{\cb}X:\cd^-\cb\ra\Tria(\cp)\cap\cd^-\ca.
\]
}
Let $\cu$ be the full subcategory of $\cd^-\cb$ formed by those $N$ such that $N\otimes^{\L}_{\cb}X\in\Tria(\cp)\cap\cd^-\ca$. It is a full triangulated subcategory of $\cd^-\cb$. Notice that, if $B=\textbf{i}P$ is the object of $\cb$ corresponding to a certain $P\in\cp$, then\[B^{\we}\otimes^{\L}_{\cb}X\cong\textbf{i}P\cong P\in\Tria(\cp)\cap\cd^-\ca.\]
This proves that $\cu$ contains the representable dg $\cb$-modules $B^{\we}$. It also proves that, since $\Tria(\cp)\cap\cd^-\ca$ is closed under small coproducts of finite extensions of objects $\Sum(\cp^+)$, then $\cu$ is closed under small coproducts of finite extensions of objects of $\Sum(\{B^{\we}\}_{B\in\cb}^+)$.  Since $\cd^-\cb$ is exhaustively generated to the left by the representable modules $B^{\we}\ko B\in\cb$, this implies that $\cu=\cd^-\cb$.
\emph{Second step: The functor $?\otimes^{\L}_{\cb}X:\cd^-\cb\ra\Tria(\cp)\cap\cd^-\ca$ is a triangle equivalence.}
To prove it we will use the techniques of \cite[Lemma 4.2]{Keller1994a}. If $B=\textbf{i}P$ is the object of $\cb$ corresponding to $P\in\cp$, we have seen already that $B^{\we}\otimes^{\L}_{\cb}X\cong P$, which is compact in $\Tria(\cp)\cap\cd^-\ca$ by hypothesis. Also, if $B=\textbf{i}P$ and $B'=\textbf{i}P'$ are objects of $\cb$, we have\begin{align}
(\cd\cb)(B^{\we},B'^{\we}[n])\arr{\sim}\H n\cb(B,B')= \nonumber \\
=(\ch\ca)(\textbf{i}P,\textbf{i}P'[n])\arr{\sim}(\cd\ca)(B^{\we}\otimes_{\cb}^{\L}X,B'^{\we}[n]\otimes_{\cb}^{\L}X). \nonumber\end{align}
Let $\cu$ be the full subcategory of $\cd^-\cb$ formed by those objects $N$ such that $?\otimes^{\L}_{\cb}X$ induces an isomorphism\[(\cd^-\cb)(B^{\we}[n],N)\arr{\sim}(\cd^-\ca)(B^{\we}[n]\otimes^{\L}_{\cb}X,N\otimes^{\L}_{\cb}X)\]
for each $B\in\cb$ and each $n\in\Z$. It is a full triangulated subcategory of $\cd^-\cb$ closed under small coproducts and containing the representable modules $B^{\we}\ko B\in\cb$. Since $\cd^-\cb$ is exhaustively generated to the left by the representable modules $B^{\we}\ko B\in\cb$, this implies that $\cu=\cd^-\cb$. Fix now an object $N\in\cd^-\cb$ and consider the full subcategory $\cv$ of $\cd^-\cb$ formed by the objects $M$ such that $?\otimes^{\L}_{\cb}X$ induces an isomorphism\[(\cd^-\cb)(M,N[n])\arr{\sim}(\cd^-\cb)(M\otimes^{\L}_{\cb}X,N[n]\otimes^{\L}_{\cb}X)\]
for each $N\in\Z$. Again, it is a full triangulated subcategory of $\cd^-\cb$ containing the representable modules and closed under small coproducts, which implies that $\cv=\cd\cb$. Therefore, we have already proved that $?\otimes^{\L}_{\cb}X$ is fully faithful.
Finally, by hypothesis, $\Tria(\cp)\cap\cd^-\ca$ is exhaustively generated to the left by the objects of $\cp$. Since they are in the essential image of the functor $?\otimes^{\L}_{\cb}X$, we deduce that it is essentially surjective.
\emph{Third step: Thanks to the second step, 1) holds if and only if the functor $?\otimes^{\L}_{\cb}X:\cd^-\cb\ra\cd^-\ca$ has a right adjoint. Let us prove that this happens if and only if $\cp$ is dually right bounded.}
The `if' part is clear. Conversely, let $G:\cd^-\ca\ra\cd^-\cb$ be a right adjoint to $?\otimes^{\L}_{\cb}X$. For simplicity, put $\R\ch om_{\ca}(X,?)=H_{X}$. Consider the diagram\[\xymatrix{& \cd^-\cb\ar@{^(->}[d]^{\iota_{\cb}}\ar@<1ex>[rr]^{?\otimes^{\L}_{\cb}X} && \cd^-\ca\ar@{^(->}[d]^{\iota_{\ca}}\ar@<1ex>[ll]^{G} \\\cd^{\leq n}\cb\ar@{^(->}[ur]\ar@<1ex>[r]^{\iota} & \cd\cb\ar@<1ex>[l]^{\tau^{\leq n}}\ar@<1ex>[rr]^{?\otimes^{\L}_{\cb}X} && \cd\ca\ar@<1ex>[ll]^{H_{X}}}
\]
where $n$ is any integer and $\iota$ is the inclusion functor. We have that\[\tau^{\leq n}\circ \iota_{\cb}\circ G\cong\tau^{\leq n}\circ H_{X}\circ\iota_{\ca}\]
since these two compositions are right adjoint to $?\otimes^{\L}_{\cb}X\circ\iota$. Let the $M$ be an object of $\cd^-\ca$ and fix an integer $n$ such that $GM\in\cd^{\leq n}\cb$. Then, we get\[\tau^{\leq n}H_{X}(M)\cong\tau^{\leq n}G(M)\cong\tau^{\leq n+i}G(M)\cong\tau^{\leq n+i}H_{X}(M)\]
for each $i\geq 0$. This implies that $\tau^{>n}(H_{X}M)\in\cd^{>n+i}\cb$ for each $i\geq 0$. In particular,\[\H j(\tau^{>n}H_{X}(M))=0
\]
for every $j\in\Z$, that is to say, $\tau^{>n}(H_{X}(M))=0$. Thus, $H_{X}(M)\in\cd^{\leq n}\cb$.\end{proof}

\subsection{Recollement of general right bounded derived categories}\label{Recollement of general right bounded derived categories}
\begin{theorem}\label{parametrization right bounded recollements}
Let $\ca$ be a dg category. The following assertions are equivalent:
\begin{enumerate}[1)]
\item $\cd^-\ca$ is a recollement of $\cd^-\cb$ and $\cd^-\cc$, for certain dg categories $\cb$ and $\cc$.\item There exist sets $\cp\ko \cq$ in $\cd^-\ca$ such that:
\begin{enumerate}[2.1)]
\item $\cp$ and $\cq$ are right bounded.
\item $\cp$ and $\cq$ are dually right bounded.
\item $\Tria(\cp)\cap\cd^-\ca$ is exhaustively generated to the left by $\cp$ and the objects of $\cp$ are compact in $\cd\ca$.\item $\Tria(\cq)\cap\cd^-\ca$ is exhaustively generated to the left by $\cq$ and the objects of $\cq$ are compact in $\Tria(\cq)\cap\cd^-\ca$.\item $(\cd\ca)(P[i],Q)=0$ for each $P\in\cp\ko Q\in\cq$ and $i\in\Z$.
\item $\cp\cup\cq$ generates $\cd\ca$.
\end{enumerate}
\end{enumerate}
\end{theorem}
\begin{proof}
$1)\Rightarrow 2)$ Consider the d\'{e}collement
\[ \xymatrix{ \cd^-\cb\ar[r]^{i_{*}=i_{!}} & \cd^-\ca\ar@/_1.2pc/[l]_{i^{*}}\ar@/_-1pc/[l]^{i^{!}}\ar[r]^{j^*=j^!} &\cd^-\cc\ar@/_1.2pc/[l]_{j_{*}}\ar@/_-1pc/[l]^{j_{!}},
}
\]
and let $(\cx',\cy',\cz')$ be the corresponding TTF triple on $\cd^-\ca$. Let $\cp$ be the set formed by all the objects $j_{!}(C^{\we})\ko C\in\cc$, and let $\cq$ be the set formed by all the objects $i_{*}(B^{\we})\ko B\in\cb$.
2.1) Notice that the coproduct $\coprod_{C\in\cc}C^{\we}$ lives in $\cd^-\cc$ and, since $j_{!}:\cd^-\cc\arr{\sim}\cx'$ is a triangle equivalence, then there exists in $\cd^-\ca$ the coproduct $\coprod_{P\in\cp}P$. Now, the claim in the proof of Lemma \ref{preservation of coproducts and compacity} implies that $\cp$ is right bounded. Similarly for $\cq$.
2.3) Since $j_{!}:\cd^-\cc\arr{\sim}\cx'$ is a triangle equivalence, then $\cx'$ is exhaustively generated to the left by the set $\cp$, whose objects are compact in $\cx'$. Then Proposition \ref{easy restriction} says that $(\cx',\cy',\cz')$ is the restriction of a TTF triple $(\cx,\cy,\cz)$ on $\cd\ca$. Moreover, $\cx=\Tria(\cp)$ and so $\cx'=\Tria(\cp)\cap\cd^-\ca$. This proves that $\Tria(\cp)\cap\cd^-\ca$ is exhaustively generated to the left by $\cp$. By using that $\cx'$ is an aisle in $\cd^-\ca$ and that $\cy'$ is closed under small coproducts in $\cd^-\ca$, we can prove that the objects of $\cp$ are compact in $\cd^-\ca$. Finally, Lemma \ref{preservation of coproducts and compacity} implies that they are also compact in $\cd\ca$.
2.4) Since $i_{*}:\cd^-\cb\arr{\sim}\cy'$ is a triangle equivalence, then $\cy'$ is exhaustively generated to the left by the set $\cq$, whose objects are compact in $\cy'$. From the proof of 2.3) we know that\[\cy'=\Tria(\cp)^{\bot}\cap\cd^-\ca.
\]
Of course, $\cq$ is contained in $\cy'$ and so $\Tria(\cq)\cap\cd^-\ca$ is contained in $\cy'$. Notice that $\Tria(\cq)\cap\cd^-\ca$ is a full triangulated subcategory of $\cy'$ containing $\cq$ and closed under small coproducts of objects of $\bigcup_{n\geq 0}\Sum(\cq^+)^{*n}$. Since $\cy'$ is exhaustively generated to the left by the set $\cq$, this implies that $\Tria(\cq)\cap\cd^-\ca=\cy'$.
2.2) From the proof of 2.3), we know that $\Tria(\cp)\cap\cd^-\ca$ is an aisle in $\cd^-\ca$. Then Proposition \ref{Morita theory for right bounded} implies that $\cp$ is dually right bounded. Similarly for $\cq$.
2.5) and 2.6) follow from the fact that $(\Tria(\cp),\Tria(\cq))$ is a t-structure on $\cd\ca$.
$2)\Rightarrow 1)$ Since the objects of $\cp$ are compact in $\cd\ca$, then Brown representability theorem implies that $(\cx,\cy):=(\Tria(\cp),\Tria(\cp)^{\bot})$ is a t-structure on $\cd\ca$. Notice that $\cy$ is closed under small coproducts and that $\tau^{\cy}$ takes a set of compact generators of $\cd\ca$ to a set of compact generators of $\cy$. Then, $\cy$ is a compactly generated triangulated category and Brown representability theorem implies that it is an aisle. Therefore,\[(\Tria(\cp),\Tria(\cp)^{\bot},(\Tria(\cp)^{\bot})^{\bot})
\]
is a TTF triple on $\cd\ca$. From conditions 2.5) and 2.6) we deduce that $\cq$ generates $\Tria(\cp)^{\bot}$. Moreover, since $\Tria(\cp)^{\bot}$ is closed under small coproducts, then $\Tria(\cq)$ is contained in $\Tria(\cp)^{\bot}$. It is an excersise to prove that the fact that $\Tria(\cq)$ is an aisle in $\cd\ca$(\cf \cite[Corollary 4.6.10]{NicolasTesis}, \cite[Corollary 3.2]{Nicolas2007}, \cite[Corollary 3.12]{Porta2007}) implies that $\Tria(\cp)^{\bot}=\Tria(\cq)$. Proposition \ref{Morita theory for right bounded} tells us that $\Tria(\cp)\cap\cd^-\ca$ and $\Tria(\cq)\cap\cd^-\ca$ are aisles in $\cd^-\ca$. Given $M\in\cd^-\ca$, consider the triangle\[M'\ra M\ra M''\ra M'[1]
\]
in $\cd^-\ca$ with $M'\in\Tria(\cp)\cap\cd^-\ca$ and $M''\in(\Tria(\cp)\cap\cd^-\ca)^{\bot}$. In particular, $M'\in\Tria(\cp)$ and $M''\in\Tria(\cp)^{\bot}=\Tria(\cq)$. This proves that\[(\Tria(\cp)\cap\cd^-\ca,\Tria(\cq)\cap\cd^-\ca)
\]
is a t-structure on $\cd^-\ca$. Similarly,
\[(\Tria(\cq)\cap\cd^-\ca,\Tria(\cq)^{\bot}\cap\cd^-\ca)
\]
is a t-structure on $\cd^-\ca$. These t-structures together form a TTF triple $(\cx',\cy',\cz')$ on $\cd^-\ca$. Finally, Proposition \ref{Morita theory for right bounded} implies that $\cx'\cong\cd^-\cc$ (for a fibrant replacement $\cc$ of $\cp$) and $\cy'\cong\cd^-\cb$ (for a fibrant replacement $\cb$ of $\cq$).\end{proof}

We will prove in Corollary \ref{parametrization right bounded hn recollements} below that conditions of assertion 2 in Theorem \ref{parametrization right bounded recollements} can be weakened under certain extra hypotheses. But first we need some preliminary results. The following one is a `right bounded' version of the proof of B.~Keller's theorem \cite[Theorem 5.2]{Keller1994a}:
\begin{proposition}\label{g implies e to the left}
Let $\cp$ be a set of objects of a triangulated category $\cd$ such that
\begin{enumerate}[1)]
\item the objects of $\cp$ are compact in $\Tria(\cp)$,
\item $\cd(P,P'[i])=0$ for each $P\ko P'\in\cp$ and $i\geq 1$,
\item small coproducts of finite extensions of objects of $\Sum(\cp^+)$ exist in $\cd$,\item for each $M\in\cd$ there exists $k_{M}\in\Z$ such that $\cd(P[n],M)=0$ for all $n<k_{M}$ and $P\in\cp$.\end{enumerate}
Then $\Tria(\cp)$ is an aisle in $\cd$ exhaustively generated to the left by $\cp$. In particular, if $\cp$ generates $\cd$, then $\Tria(\cp)=\cd$.\end{proposition}
\begin{proof}
We include the proof for the sake of completeness. Let $M\in\cd$. We know that if $\cd(P[n],M)\neq 0$ for some $P\in\cp$, then $n\geq k_{M}$. Since $\cp$ is a set, there exists an object\[P_{0}\in\Sum(\cp^+[k_{M}])
\]
and a morphism $\pi_{0}:P_{0}\ra M$ inducing a surjection
\[\pi_{0}^{\we}:\cd(P[n],P_{0})\ra\cd(P[n],M)
\]
for each $P\in\cp\ko n\in\Z$. Indeed, one can take
\[P_{0}:=\coprod_{P\in\cp\ko n\geq k_{M}}P[n]^{(\cd(P[n],M))}.
\]
Now, we will inductively construct a commutative diagram
\[\xymatrix{P_{0}\ar[r]^{f_{0}}\ar[dr]_{\pi_{0}} & P_{1}\ar[r]^{f_{1}}\ar[d]^{\pi_{1}} & \dots\ar[r] & P_{q}\ar[r]^{f_{q}}\ar[dll]^{\pi_{q}} & \dots\ko q\geq 0 \\& M &&&
}
\]
such that:
\begin{enumerate}[a)]
\item $P_{q}\in\Sum(\cp^+[k_{M}])^{*q}$,
\item $\pi_{q}$ induces a surjection
\[\pi_{q}^{\we}:\cd(P[n],P_{q})\ra\cd(P[n],M)
\]
for each $P\in\cp\ko n\in\Z$.
\end{enumerate}
Suppose for some $q\geq 0$ we have constructed $P_{q}$ and $\pi_{q}$. Consider the triangle\[C_{q}\arr{\alpha_{q}}P_{q}\arr{\pi_{q}}M\ra C_{q}[1]
\]
induced by $\pi_{q}$. By applying $\cd(P[n],?)$ we get a long exact sequence\[\dots\ra\cd(P[n+1],M)\ra\cd(P[n],C_{q})\ra\cd(P[n],P_{q})\ra\dots
\]
If $\cd(P[n],C_{q})\neq 0$, then either $\cd(P[n+1],M)\neq 0$ or $\cd(P[n],P_{q})\neq 0$. In the first case, we would have $n\geq k_{M}-1$. In the second case we would have $\cd(P[n],P'[m])\neq 0$ for some $P'\in\cp\ko m\geq k_{M}$, and so $n\geq m\geq k_{M}$. Therefore, $\cd(P[n],C_{q})\neq 0$ implies $n\geq k_{M}-1$. This allows us to take\[Z_{q}\in\Sum(\cp^+[k_{M}-1])
\]
together with a morphism $\beta_{q}:Z_{q}\ra C_{q}$ inducing a surjection\[\beta_{q}^{\we}:\cd(P[n],Z_{q})\ra \cd(P[n],C_{q})
\]
for each $P\in\cp\ko n\in\Z$. Define $f_{q}$ by the triangle
\[Z_{q}\arr{\alpha_{q}\beta_{q}}P_{q}\arr{f_{q}}P_{q+1}\ra Z_{q}[1]
\]
Since $\pi_{q}\alpha_{q}=0$, there exists $\pi_{q+1}:P_{q+1}\ra M$ such that $\pi_{q+1}f_{q}=\pi_{q}$. Notice that, since\[Z_{q}[1]\in\Sum(\cp^+[k_{M}]),
\]
then
\[ P_{q+1}\in\Sum(\cp^+[k_{M}])^{*(q+1)}.
\]
Also, the surjectivity required for $\pi^{\we}_{q+1}$ follows from the surjectivity guaranteed for $\pi^{\we}_{q}$. Define $P_{\infty}$ to be the Milnor colimit of the sequence $f_{q}\ko q\geq 0$:\[\coprod_{q\geq 0}P_{q}\arr{\phi}\coprod_{q\geq 0}P_{q}\arr{\psi} P_{\infty}\ra \coprod_{q\geq 0}P_{q}[1].\]
Consider the morphism
\[\theta=\scriptsize{\left[\begin{array}{ccc}\pi_{0}&\pi_{1}&\dots\end{array}\right]}:\coprod_{q\geq 0}P_{q}\ra M.\]
Since $\pi_{q+1}f_{q}=\pi_{q}$ for every $q\geq 0$, we have $\theta\phi=0$, which induces a morphism $\pi_{\infty}:P_{\infty}\ra M$ such that $\pi_{\infty}\psi=\theta$. If we prove that $\pi_{\infty}$ induces an isomorphism\[\pi^{\we}_{\infty}:\cd(P[n],P_{\infty})\arr{\sim}\cd(P[n],M)
\]
for every $P\in\cp\ko n\in\Z$, then we would have
\[\cd(P[n],\cone(\pi_{\infty}))=0
\]
for every $P\in\cp\ko n\in\Z$, that is to say
\[\cone(\pi_{\infty})\in\Tria(\cp)^{\bot}.
\]
Therefore, we would have proved that $\Tria(\cp)$ is an aisle in $\cd$. Also, if $M\in\Tria(\cp)$, in the triangle\[P_{\infty}\arr{\pi_{\infty}}M\ra\cone(\pi_{\infty})\ra P_{\infty}[1]
\]
we would have that $P_{\infty}\ko M\in\Tria(\cp)$, which implies
\[\cone(\pi_{\infty})\in\Tria(\cp).
\]
Therefore, $\cone(\pi_{\infty})=0$ and so $\pi_{\infty}$ is an isomorphism. Thus, we would have proved that for every object of $\Tria(\cp)$ there exists an integer $k_{M}$ and a triangle\[\coprod_{q\geq 0}P_{q}\ra\coprod_{q\geq 0}P_{q}\ra M[-k_{M}]\ra\coprod_{q\geq 0}P_{q}[1]\]
with $P_{q}\in\Sum(\cp^+)^{*q}\ko q\geq 0$. In particular, we would have that $\Tria(\cp)$ is exhaustively generated to the left by $\cp$.
Let us prove the bijectivity of $\pi^{\we}_{\infty}$. The surjectivity follows from the identity $\pi^{\we}_{\infty}\psi^{\we}=\theta^{\we}$ and the fact that $\theta^{\we}$ is surjective (thanks to the surjectivity of the $\pi_{q}^{\we}\ko q\geq 0$ and the  compactness of the $P\in\cp$). Now consider the commutative diagram\[\xymatrix{\coprod_{q\geq 0}\cd(P[n],P_{q})\ar[r]^{\phi^{\we}} & \coprod_{q\geq 0}\cd(P[n],P_{q})\ar[r]^{\psi^{\we}}\ar[dr]_{\theta^{\we}} & \cd(P[n],P_{\infty})\ar[r]\ar[d]^{\pi^{\we}_{\infty}} & 0 \\&& \cd(P[n],M) &
}
\]
The map $\psi^{\we}$ is surjective since the map
\[\phi[1]^{\we}:\coprod_{q\geq 0}\cd(P[n],P_{q}[1])\ra \coprod_{q\geq 0}\cd(P[n],P_{q}[1])\]
is injective. If we prove that the kernel of $\theta^{\we}$ is contained in the image of $\phi^{\we}$, then we would have the injectivity of $\pi^{\we}_{\infty}$ by an easy diagram chase. Let\[g=\left[\begin{array}{cccccc}g_{0} & g_{1} & \dots & g_{s}&0&\dots\end{array}\right]^{t}:P[n]\ra\coprod_{q\geq 0}P_{q}\]
be an element of the kernel of $\theta^{\we}$. Then
\[\pi_{0}g_{0}+\dots +\pi_{s}g_{s}=0
\]
implies
\[\pi_{s}(f_{s-1}\dots f_{0}g_{0}+ f_{s-1}\dots f_{1}g_{1}+\dots +g_{s})=0\]
and so the morphism
\[f_{s-1}\dots f_{0}g_{0}+ f_{s-1}\dots f_{1}g_{1}+\dots +g_{s}
\]
factors through $\alpha_{s}$:
\[f_{s-1}\dots f_{0}g_{0}+ f_{s-1}\dots f_{1}g_{1}+\dots +g_{s}=\alpha_{s}\gamma_{s}:P[n]\ra C_{s}\ra P_{s}.\]
By construction of $Z_{s}$ we have that $\gamma_{s}$ factors through $\beta_{s}$, and so\[f_{s-1}\dots f_{0}g_{0}+ f_{s-1}\dots f_{1}g_{1}+\dots +g_{s}=\alpha_{s}\beta_{s}\xi_{s}.\]
This implies
\[f_{s}\dots f_{0}g_{0}+ f_{s}\dots f_{1}g_{1}+\dots +f_{s}g_{s}=f_{s}\alpha_{s}\beta_{s}\xi_{s}=0,\]
since $f_{s}\alpha_{s}\beta_{s}=0$ by construction of $f_{s}$. Therefore, the morphism\[h:P[n]\ra\coprod_{q\geq 0}P_{q}
\]
with non-vanishing components
\[P[n]\ra P_{r}\ra\coprod_{q\geq 0}P_{q}
\]
induced by
\[g_{r}+\dots + f_{r-1}\dots f_{1}g_{1}+f_{r-1}\dots f_{0}g_{0}: P[n]\ra P_{r}\]
with $0\leq r\leq s$, satisfies $\varphi^{\we}(h)=g$.
\end{proof}

\begin{corollary}\label{right bounded compactly generated implies exhaustive}Let $\ca$ be a dg category and let $\cp$ be a set of objects of $\cd^-\ca$ such that:\begin{enumerate}[1)]
\item it is both right bounded and dually right bounded,
\item its objects are compact in $\Tria(\cp)\cap\cd^-\ca$,
\item $(\cd\ca)(P,P'[i])=0$ for each $P\ko P'\in\cp$ and $i\geq 1$.
\end{enumerate}
Then $\Tria(\cp)\cap\cd^-\ca$ is exhaustively generated to the left by $\cp$.\end{corollary}
\begin{proof}
Put $\cd:=\Tria(\cp)\cap\cd^-\ca$. Since $\cp$ is right bounded, then $\cp$ is contained in $\cd$ and for each integer $k$ small coproducts of finite extensions of $\Sum(\cp^+[k])$ are in $\cd$. Also, Proposition \ref{characterization of dually right bounded} guarantees that for each $M\in\cd$ there exists an integer $k_{M}$ such that $\cd(P[n],M)=0$ for each $P\in\cp$ and $n<k_{M}$. Therefore, we can apply Proposition \ref{g implies e to the left}.\end{proof}

\begin{corollary}\label{parametrization right bounded hn recollements}
Let $\ca$ be a dg category. The following assertions are equivalent:
\begin{enumerate}[1)]
\item $\cd^-\ca$ is a recollement of $\cd^-\cb$ and $\cd^-\cc$, for certain dg categories $\cb$ and $\cc$ with cohomology concentrated in non-positive degrees.\item There exist sets $\cp\ko \cq$ in $\cd^-\ca$ such that:
\begin{enumerate}[2.1)]
\item $\cp$ and $\cq$ are right bounded.
\item $\cp$ and $\cq$ are dually right bounded.
\item The objects of $\cp$ are compact in $\cd\ca$ and satisfy
\[(\cd\ca)(P,P'[i])=0
\]
for all $P\ko P'\in\cp$ and $i\geq 1$.
\item The objects of $\cq$ are compact in $\Tria(\cq)\cap\cd^-\ca$ and satisfy\[(\cd\ca)(Q,Q'[i])=0
\]
for all $Q\ko Q'\in\cq$ and $i\geq 1$.
\item $(\cd\ca)(P[i],Q)=0$ for each $P\in\cp\ko Q\in\cq$ and $i\in\Z$.
\item $\cp\cup\cq$ generates $\cd\ca$.
\end{enumerate}
\end{enumerate}
\end{corollary}
\begin{proof}
$1)\Rightarrow 2)$ Is similar to the corresponding implication in Theorem \ref{parametrization right bounded recollements}. The fact that the dg categories $\cb$ and $\cc$ have cohomology concentrated in non-positive degrees is reflected in the fact that\[(\cd\ca)(P,P'[i])=(\cd\ca)(Q,Q'[i])=0
\]
for each $P\ko P'\in\cp\ko Q\ko Q'\in\cq$ and $i\geq 1$.

$2)\Rightarrow 1)$ Thanks to Corollary \ref{right bounded compactly generated implies exhaustive}, conditions 2.3 and 2.4 of Theorem \ref{parametrization right bounded recollements} are satisfied. Therefore, that Proposition (and its proof) ensures that $\cd^-\ca$ is a recollement of $\cd^-\cb$ and $\cd^-\cc$, where $\cb$ is a fibrant replacement of $\cq$ and $\cc$ is a fibrant replacement of $\cp$. Finally, the fact that\[(\cd\ca)(P,P'[i])=(\cd\ca)(Q,Q'[i])=0
\]
for each $P\ko P'\in\cp\ko Q\ko Q'\in\cq$ and $i\geq 1$ implies that $\cb$ and $\cc$ have cohomology concentrated in non-positive degrees.\end{proof}

\subsection{Recollement of right bounded derived categories of algebras}\label{Recollement of right bounded derived categories of algebras}
\begin{definition}
Let $A$ be an ordinary algebra. If $M$ is a complex of $A$-modules, the \emph{graded support} of $M$ is the set ofintegers $i\in\Z$ such that $M^i\neq 0$. In case $M$ is a bounded complex, we consider\[w(M):=sup\{i\in\Z\mid M^i\neq 0\}-inf\{i\in\Z\mid M^i\neq 0\}+1
\]
and call it the \emph{width} of $M$. Suppose now that $P$ is a bounded complex of projective $A$-modules, so that $P$ is a dually right bounded object of $\cd^-A$ (\cf Proposition \ref{characterization of dually right bounded}),  and $M\in\cd^-A$ is any object of the right bounded derived category. Unless $M\in\Tria_{\cd A}(P)^\perp$, there is a well-defined integer $k_M:=inf\{n\in\Z\mid (\cd A)(P[n],M)\neq 0\}$.\end{definition}

\begin{lemma} \label{route to Konig2}
Let $A$ be an ordinary algebra. Let $P$ be a bounded complex of projective $A$-modules such that $(\cd A)(P,P[i])=0$, for all $i>0$, and the canonicalmorphism $(\cd A)(P,P[i])^{(\Lambda)}\ra(\cd A)(P,P[i]^{(\Lambda )})$ is an isomorphism, for every integer $i$ and every set $\Lambda$.Let $M$ be an object of $\Tria_{\cd A}(P)\cap\cd^-A$. There exists a sequence of inflations $0=P_{-1}\ra P_0\ra P_1\ra ...$ in $\cc A$, whose colimit is denoted by $P_{\infty}$, satisfying the following properties:\begin{enumerate}[1)]
\item $P_\infty$ is isomorphic to $M$ in $\cd A$.
\item $P_n/P_{n-1}$ belongs to $\Sum(\{P\}^+[k_{M}+n])$, for each $n\geq 0$.\item If $n\geq w(P)-k_M$ the graded supports of $P$ and $P_\infty/P_n$ are disjoint.\end{enumerate}
\end{lemma}
\begin{proof}
Imitating the proof of Proposition \ref{g implies e to the left}, we shall construct a filtration satisfying conditions 2) and 3), leaving for the last moment the verification of condition 1).
\emph{First step: condition 2).} Note that in the proof of that proposition, we start with $P_0\in\Sum(P[i]:$ $i\geq k_M)$ and then, at each step, $P_{q+1}$ appears in a triangle\[Z_q\arr{\alpha_q\beta_q}P_q\arr{f_q}P_{q+1}\ra Z_q[1],
\]
where $Z_q$ is a coproduct of shifts $P[i]$, with $i\geq k_{M}-1$. Working in $\cc A$ and bearing in mind that $Z_q$ is cofibrant (it is a right bounded complex of projective $A$-modules), we can assume without loss of generality that $f_q$ is the mapping cone of a cochain map $Z_q\ra P_q$ and, as a consequence, that $f_q$ is an inflation in $\cc A$ appearing in a conflation\[P_q\arr{f_q}P_{q+1}\ra Z_q[1],
\]
where $Z_q[1]$ is a coproduct in $\cc A$ of shifts $P[i]$, $i\geq k_M$. We shall prove by induction on $q\geq 0$ that one can choose$Z_q[1]\in\Sum(\{P\}^+[q+1+k_{M}])$ or, equivalently, that $Z_q\in\Sum(\{P\}^+[q+k_{M}])$. Since $Z_q$ is defined via a map $\beta_{q}:Z_{q}\ra C_{q}$ such that\[\beta^{\we}_{q}:(\cd A)(P[i],Z_{q})\ra(\cd A)(P[i],C_{q})
\]
is surjective for all $i\in\Z$, our task reduces to prove that $(\cd A)(P[i],C_q)\neq 0$ implies $i\geq q+k_M$. We leave as an exercise checking that for $q=0$. Provided it is true for $q-1$, we apply the homological functor $(\cd A)(P[i],?)$ to the triangle\[Z_{q-1}\arr{\beta_{q-1}}C_{q-1}\arr{u_{q-1}}C_q\ra Z_{q-1}[1]
\]
and, bearing in mind that $(\cd A)(P[i],Z_{q-1})\ra(\cd A)(P[i],C_{q-1})$ is surjective, we get that $(\cd A)(P[i],C_{q})\ra(\cd A)(P[i],Z_{q-1}[1])$ is injective. As a consequence, the inequality $(\cd A)(P[i],C_{q})\neq 0$ implies that $(\cd A)(P[i],Z_{q-1}[1])\neq 0$ and the induction hypothesis guarantees that $Z_{q-1}$ is a coproduct of shifts $P[j]$, with $j\geq q-1+k_M$. Then $(\cd A)(P[i],C_{q})\neq 0$ implies that $0\neq(\cd A)(P[i],P[j+1])=(\cd A)(P,P[j+1-i])$,for some $j\geq q-1+k_M$. Then $i\geq q+k_M$ as desired. In conclusion, we can view the map $f_q:P_q\ra P_{q+1}$ as an inflation in $\cc A$ whose cokernel is isomorphic in $\cc A$ to a coproduct of shifts $P[i]$, with $i\geq q+1+k_M$.

\emph{Second step: condition 3).} If now $n\geq 0$ is any natural number, then $P_\infty/P_n$ admits a filtration\[0=P_n/P_n\ra P_{n+1}/P_n\ra ...
\]
in $\cc A$, where the quotient of two consecutive factors is a coproduct of shifts $P[i]$, with $i\geq n+k_M$. If $n\geq w(P)-k_M$, then any such index $i$ satisfies $i\geq w(P)$ and then the graded supports of $P$ and $P[i]$ are disjoint. As a result the graded supports of $P$ and $P_\infty/P_n$ are disjointwhenever $n\geq w(P)-k_M$.

\emph{Third step: condition 1).} Finally, in order to prove condition 1), notice that the argument in the final part of the proof of Proposition \ref{g implies e to the left} can be repeated, as soon as we are able to prove that the canonical morphism $\coprod_{n\geq 0}(\cd A)(P[i],P_n)\ra(\cd A)(P[i],\coprod_{n\geq 0}P_n)$ is an isomorphism, for every integer $i\in\Z$. It is not difficult to reduce that to the case in which $i=0$. For that we fix $n\geq w(P)-k_M$ large enough so that also the graded supports of $P[1]$ and $P_\infty /P_n$ are disjoint. Then we get a conflation in $\cc A$\[\left(\coprod_{k\leq n}P_k\right)\oplus\left(\coprod_{k>n} P_n \right)\ra\coprod_{k\geq 0}P_k\ra\coprod_{k>n}P_k/P_n.\]
That conflation of $\cc A$ gives rise to the corresponding triangle of $\cd A$. But the right term in the above conflation has a graded support which isdisjoint with those of $P$ and  $P[1]$. That implies that
\[(\cd A)(P,\coprod_{k>n}P_k/P_n)=0=(\cd A)(P,\coprod_{k>n}P_k/P_n[-1])\]
and also
\[\coprod_{k>n}(\cd A)(P,P_k/P_n)=0=\coprod_{k>n}(\cd A)(P,P_k/P_n[-1]).\]
We then get a commutative diagram with horizontal isomorphisms:
\[\xymatrix{\left(\coprod_{k\leq n}(\cd A)(P,P_{k})\right)\oplus\left(\coprod_{k>n}(\cd A)(P,P_{n})\right)\ar[d]^{\can}\ar[rr]^{\ \ \ \ \ \ \ \ \ \sim} && \coprod_{k\geq 0}(\cd A)(P,P_{k})\ar[d]^{\can} \\\left((\cd A)(P,\coprod_{k\leq n}P_{k})\right)\oplus\left((\cd A)(P,\coprod_{k>n}P_{n})\right)\ar[rr]^{\ \ \ \ \ \ \ \ \ \sim} && (\cd A)(P,\coprod_{k\geq 0}P_{k})}
\]
The proof will be finished if we are able to prove, for any fixed natural number $n$, that $(\cd A)(P[i],?)$ preserves small coproducts of objects in $\Sum(\{P\}^+)^{*n}$ for every $i\in\Z$. Let us prove it. From the hypotheses on $P$ and the fact if $i>w(P)$ then $(\cd A)(P,P[i]^{(\Lambda )})=0$ for every set $\Lambda$, one readily sees that, for every integer $m$ and every family of exponent sets $(\Lambda_i)_{i\geq m}$, the canonical morphism$\coprod_{i\geq m}(\cd A)(P,P[i])^{(\Lambda_i)}\ra(\cd A)(P,\coprod_{i\geq m}P[i]^{(\Lambda_i)})$ is an isomorphism. Our goal is then attained for $n=0$ and an easy induction argument gets the job done for every $n\geq 0$.\end{proof}

\begin{definition}\label{exceptional}
An object $M$ of a (tipically compactly generated) triangulated category $\cd$ is \emph{exceptional} if $\cd(M,M[i])=0$ for every integer $i\neq 0$.\end{definition}

Now we can deduce the following theorem:

\begin{theorem} \label{ourKonig}
Let $A$, $B$ and $C$ be ordinary algebras. The following assertions are equivalent:\begin{enumerate}[1)]
\item $\cd^-A$ is a recollement of $\cd^-C$ and $\cd^-B$.
\item There are two objects $P\ko Q\in\cd^-A$ satisfying the following properties:\begin{enumerate}[2.1)]
\item There are isomorphisms of algebras $C\cong(\cd A)(P,P)$ and $B\cong(\cd A)(Q,Q)$.\item $P$ is exceptional and isomorphic in $\cd A$ to a bounded complex of finitely generated projective $A$-modules.\item For every set $\Lambda$ and every non-zero integer $i$ we have $(\cd A)(Q,Q^{(\Lambda )}[i])=0$, the canonical map $(\cd A)(Q,Q)^{(\Lambda)}\ra(\cd A)(Q,Q^{(\Lambda )})$ is an isomorphism, and $Q$ is isomorphic in $\cd A$ to a bounded complex of projective $A$-modules.\item $(\cd A)(P,Q[i])=0$ for all $i\in\Z$.
\item $P\oplus Q$ generates $\cd A$.
\end{enumerate}
\end{enumerate}
\end{theorem}
\begin{proof}
$1)\Rightarrow 2)$ is a particular case of the proof of the corresponding implication in Corollary \ref{parametrization right bounded hn recollements}, where we take into account Corollary \ref{dually right bounded are bounded complexes of projectives} and the additional consideration that the dg categories are in this case ordinary algebras, whence having cohomology concentrated in degree zero.
$2)\Rightarrow 1)$ Taking $\cp=\{P\}$ and $\cq=\{Q\}$, one readily sees that these one-point sets satisfy conditions 2.1, 2.2, 2.3, 2.5 and 2.6 of Corollary \ref{parametrization right bounded hn recollements}. As for condition 2.4 it only remains to prove that $Q$ is compact in $\Tria_{\cd A}(Q)\cap\cd^-A$. For this,let $(M_{j})_{j\in J}$ be a family of objects in $\Tria(Q)\cap\cd^-A$ having a coproduct, say $M$,  in that subcategory and denote by $q_j:M_{j}\ra M$ the injections. Of course, we have that $sup\{i\in\Z\mid \H i(M_{j})\neq 0\}\leq sup\{i\in\Z\mid \H i(M)\neq 0\}$, for every $j\in J$. Thenthe coproduct $\coprod_{j\in J} M_{j}$ of the family in $\cd A$ belongs to $\cd^-A$ and thus to $\Tria(Q)\cap\cd^-A$. This easily implies that$M\cong\coprod_{j\in J}M_{j}$ and the injection $q_j:M_{j}\ra M$ gets identified with the canonical injection $M_{j}\ra\coprod_{k\in J}M_{k}$.For each $j\in J$ we consider the complex $Q_{j,\infty}$ and the filtration\[0=Q_{j,-1}\ra Q_{j,0}\ra Q_{j,1}\ra ...
\]
given by Lemma \ref{route to Konig2} for $M_{j}$, where we have replaced the letter ``P'' by the letter ``Q'' to avoid confusion with theobject $P$. Notice that $k_{M}\leq k_{M_{j}}$ for every $j\in J$. Therefore, the integer $r:=inf\{k_{M_{j}}\}_{j\in J}$ is well defined. If we fix $n\in\mathbf{N}$ such that $n+r>w(Q)$, then $n>w(Q)-k_{M_{j}}$. Notice that \cite[Lemma 5.3]{Keller1990} implies that a countable composition of inflations of $\cc A$ is again an inflation of $\cc A$. Then, for every $j\in J$ we get a conflation in $\cc A$,\[Q_{j,n}\ra Q_{j,\infty}\ra Q_{j,\infty}/Q_{j,n},
\]
By Lemma \ref{route to Konig2}, the right term of this conflation has a graded support which is disjoint with that of $Q$ and $Q[1]$ (enlarging $n$ if necessary). Then we get a commutative diagram:\[\xymatrix{\coprod_{J}(\cd A)(Q,Q_{j,n})\ar[d]^{\can}\ar[r]^{\sim} & \coprod_{J}(\cd A)(Q,Q_{j,\infty})\ar[d]^{\can}\ar[r]^{\sim} & \coprod_{J}(\cd A)(Q,M_{j})\ar[d]^{\can} \\(\cd A)(Q,\coprod_{J}Q_{j,n})\ar[r]^{\sim} & (\cd A)(Q,\coprod_{J}Q_{j,\infty})\ar[r]^{\sim} & (\cd A)(Q,\coprod_{J}M_{j})}
\]
The fact that the leftmost vertical map is a bijection has been proved in the third step of the proof of Lemma \ref{route to Konig2}, and so we are done.\end{proof}

\section{More than an exceptional object}\label{More than an exceptional object}
\subsection{The mismatch}\label{The mismatch}

Theorem \ref{ourKonig} is very close to the following theorem of S.~K\"{o}nig \cite[Theorem 1]{Konig1991}:
\begin{thm} \label{Konig}
Let $A$, $B$ and $C$ be ordinary algebras. The following assertions are equivalent:\begin{enumerate}[1)]
\item $\cd^-A$ is a recollement of $\cd^-C$ and $\cd^-B$.
\item There are two objects $P\ko Q\in\cd^-A$ satisfying the following properties:\begin{enumerate}[2.1)]
\item There are isomorphisms of algebras $C\cong(\cd A)(P,P)$ and $B\cong(\cd A)(Q,Q)$.\item $P$ is exceptional and isomorphic in $\cd A$ to a bounded complex of finitely generated projective $A$-modules.\item $Q$ is exceptional, it is isomorphic in $\cd A$ to a bounded complex of projective $A$-modules and the functor $\Hom_{A}(Q,?):\Mod A\ra\Mod k$ preserves small coproducts of copies of $Q$.\item $(\cd A)(P,Q[i])=0$ for all $i\in\Z$.
\item $P\oplus Q$ generates $\cd A$.
\end{enumerate}
\end{enumerate}
\end{thm}

The reader will have noticed that we changed S.~K\"{o}nig's condition
that $Q$ is exceptional for the stronger condition that $(\cd
A)(Q,Q[i]^{(\Lambda )})=0$ , for all $i\neq 0$ and all sets
$\Lambda$. In what follows we will show that this stronger
condition is needed in order for the theorem to be valid.

\subsection{Some results on countable von Neumann regular algebras}

We thank J.~Trlifaj for giving us an example \cite[Lemma 6.3.14 and Example 6.3.15]{NicolasTesis} that was at the basis for the following developement.
\begin{lemma}\label{is hereditary}
If $A$ is a countable von Neumann regular algebra, then it is hereditary on both sides.\end{lemma}
\begin{proof}
Since $A$ is countable its pure global dimension on either side is smaller or equal than $1$ \cite[Th\`{e}orem 7.10]{GrusonJensen} and, since $A$ is von Neumann regular, we conclude that $A$ is hereditary on both sides \cite[Proposition 10.3]{GrusonJensen}.\end{proof}

\begin{lemma}\label{preserves coproducts}
Let $A$ be a countable simple von Neumann regular algebra which is not
semisimple. If $Q$ is an injective $A$-module then the functor
$\Hom_{A}(Q,?):\Mod A\ra\Mod k$ preserves small coproducts.
\end{lemma}
\begin{proof}
\emph{First step: the countable sequence of submodules.}
of submodules of $Q$ where $Q_{n}:=f^{-1}(\coprod_{i=0}^{n}M_{i})$. Notice that $Q=\bigcup_{n\in\N}Q_{n}$ and that for every $n\in\N$ we have $Q_{n}\neq Q$. This implies that we can choose a sequence $n_{0}<n_{1}<\dots$ of natural numbers such that $Q_{i}$ is strictly contained in $Q_{i+1}$ whenever $i=n_{t}$ for some $t\in\N$.
\emph{Second step: the countable sequence of idempotents.} Let $e_{t}\ko t\in\N$, be a sequence of mutually orthogonal non-zero idempotents of $A$. Since $Ae_{t}A$ is a non-zero two-sided ideal of the simple algebra $A$, then $Ae_{t}A=A$ and so $Q_{n_{t}}=Q_{n_{t}}A=Q_{n_{t}}e_{t}A$. Therefore, for each $t\in\N$ there exists an element $x_{t}\in Q_{n_{t}}e_{t}$ which does not belong to $Q_{n_{t-1}}$.
\[g:\bigoplus_{t\in\N}e_{t}A\ra Q\ko \sum_{t\in\N}a_{t}\mapsto\sum_{t\in\N}x_{t}a_{t}.\]
Since $Q$ is injective, $g$ extends to $A$ and so there exists an element $x\in Q$ such that for every $t\in\N$ we have that $g(e_{t})=x_{t}e_{t}=x_{t}=xe_{t}$. If $s$ is a natural number such that $x\in Q_{n_{s}}$, then $x_{t}\in Q_{n_{s}}$ for every $t\in\N$, which contradicts the choice of the elements $x_{t}$.\end{proof}

\begin{lemma}\label{sin algunas extensiones}
Let $A$ be a countable simple von Neumann regular algebra which is not right N\oe therian. If $Q$ is an injective cogenerator of $\Mod A$ containing an isomorphic copy of every cyclic module, then $\Ext^{1}_{A}(Q,Q^{(\N)})\neq 0$.\end{lemma}
\begin{proof}
\emph{First step: $Q^{(\N)}$ is not injective.}
\emph{Second step: $\Ext^{1}_{A}(Q,Q^{(\N)})\neq 0$.} Since $Q^{(\N)}$ isnot injective, Baer's criterion implies that \cite[Theorem
18.3]{AndersonFuller1992} there exists a cyclic $A$-module $M$ such that
$\Ext^{1}_{A}(M,Q^{(\N)})\neq 0$. We fix a monomorphism
$j:M\longrightarrow Q$, which we view as an inclusion. Now, by applying
$\Hom_{A}(Q,?)$ and $\Hom_{A}(M,?)$ to the minimal injective coresolution\[0\ra Q^{(\N)}\ra E(Q^{(\N)})\ra E'\ra 0
\]
we get a commutative diagram with exact rows
\[\xymatrix{\Hom_{A}(Q,E')\ar[r]\ar[d] & \Ext^{1}_{A}(Q,Q^{(\N)})\ar[d]\ar[r] & 0 \\\Hom_{A}(M,E')\ar[r] & \Ext^{1}_{A}(M,Q^{(\N)})\ar[r] & 0
}
\]
where the left vertical arrow is the restriction map, and it is surjective because $E'$ is injective. Then, the right vertical arrow is surjective, which implies that $\Ext^{1}_{A}(Q,Q^{(\N)})\neq 0$.\end{proof}

\begin{lemma}\label{is a generator}
Let $A$ be a countable simple von Neumann regular algebra which is not right N\oe therian, and let $Q$ be an injective cogenerator of $\Mod A$ containing an isomorphic copy of every cyclic module. Then an $A$-module $M$ is zero whenever\[\Hom_{A}(Q,M)=\Ext^{1}_{A}(Q,M)=0.
\]
\end{lemma}
\begin{proof}
Consider a minimal injective coresolution $0\ra M\ra E(M)\ra E'(M)\ra 0$ of an $A$-module $M$ such that $\Hom_{A}(Q,M)=\Ext^{1}_{A}(Q,M)=0$.
\emph{First step: If $M\neq 0$ then $E'(M)\neq 0$.} Indeed, if $E'(M)=0$then $M$ is injective and so it contains the injective envelope of any
non-zero cyclic submodule of $M$, which would be a non-zero direct summandof $Q$ and $M$. This implies $\Hom_{A}(Q,M)\neq 0$, which is a
contradiction.

\emph{Second step: $M=0$.} Suppose not and let $C$ be a non-zero cyclicsubmodule of $E'(M)$, so that its injective envelope $Q':=E(C)$ is a
direct summand of $E'(M)$. Fix a section $v:Q'\ra E'(M)$. Since
$\Ext^{1}_{A}(Q',M)=0$, there exists a morphism of $A$-modules $f:Q'\ra E(M)$ which fits in the following commutative diagram
\[\xymatrix{ &&& Q'\ar[d]^{v}\ar[dl]_{f}& \\
0\ar[r]& M\ar[r] & E(M)\ar[r]_{p} & E'(M)\ar[r] & 0
}
\]
Then $f$ is a monomorphism and $f(Q')$ is a direct summand of $E(M)$ isomorphic to $Q'$ and such that $p$ induces and isomorphism $\pi:f(Q')\ra v(Q')$. Hence, we can rewrite the short exact sequence above as\[\xymatrix{0\ra M\ra E\oplus f(Q')\ar[rr]^{\scriptsize{\ \ \ \ \ \ \ \ \ \ \left[\begin{array}{cc}\alpha & 0 \\ \beta & \pi\end{array}\right]}} && E'\oplus Q'\ar[r] & 0}
\]
Notice that in this short exact sequence the kernel of the epimorphism, which is isomorphic to $M$, intersects in $0$ with $0\oplus f(Q')$. This implies that $M$ is not essential in $E(M)$, which is absurd.\end{proof}

\begin{example}\label{TrlifajExample}
Any countable direct limit of countable simple Artinian algebras is a countable simple von Neumann regular algebra which is not right N\oe therian. A typical case is given as follows. Consider the direct limit $\lid \cm_{2^n\times 2^n}(\mathbb{K})$, where $\mathbb{K}$ is a countable field and the ring morphism$\cm_{2^n\times 2^n}(\mathbb{K})\ra \cm_{2^{n+1}\times 2^{n+1}}(\mathbb{K})$ maps the matrix $U$ onto the matrix given by the block decomposition$\left[\begin{array}{cc} U & 0\\ 0 & U \end{array}\right]$.
\end{example}

\subsection{A counterexample}

Let $H$ be a countable simple von Neumann regular algebra which is not
right N\oe therian. Let $\cl(H)$ be the family of right ideals of $H$ andlet $Q'$ be the injective envelope of $\coprod_{I\in\cl(H)}H/I$.

\begin{remark}
Notice that $\op{End}_{H}(Q')$ is not countable for there exist two
obvious injective maps
\[\prod_{I\in\cl(H)}\op{End}_{H}(H/I)\ra\op{End}_{H}(\coprod_{I\in\cl(H)}H/I)\ra\op{End}_{H}(Q').\]
\end{remark}

Let $C$ be any unital subalgebra of $\End_{H}(Q')$. Take $A$ to be the triangular matrix algebra\[A:=\left[\begin{array}{cc}C&Q'\\ 0&H\end{array}\right].
\]
The category $\Mod A$ admits a nice description in terms of $\Mod C$ and $\Mod H$ (\cf \cite[subsection III.2]{AuslanderReitenSmalo1995}) that we will use without explicit mention. Consider the idempotent\[e:=\left[\begin{array}{cc}1&0\\ 0&0\end{array}\right].
\]
Since $AeA=eA$, it turns out that the idempotent ideal $I:=AeA$ is a projective right $A$-module and so the canonical projection\[\pi:A\ra A/I
\]
is a homological epimorphism (\cf \cite[Example 5.3.4]{NicolasTesis}). Of course, we have isomorphisms of algebras $A/I\cong H$ and $\End_{H}(I)\cong eAe\cong C$.Therefore, $\cd A$ is a recollement of $\cd H$ and $\cd C$ as follows:
\[\xymatrix{\cd H\ar[rr] && \cd A\ar@/_1pc/[ll]^{}\ar@/_-1pc/[ll]^{}\ar[rr] && \cd C.\ar@/_1pc/[ll]\ar@/_-1pc/[ll]}
\]
It induces a TTF triple $(\cx,\cy,\cz)$ on $\cd A$ where $\cx=\Tria_{\cd A}(I)$ and $\cy$ consists of those complexes isomorphic in $\cd A$ to complexes of $H$-modulesregarded as $A$-modules. Moreover, Example \ref{restriction in examples} tells us that this TTF triple restricts to a TTF triple on $\cd^-A$ which expresses $\cd^-A$ as a recollement of $\cd^-H$ and $\cd^-C$:\[\xymatrix{\cd^- H\ar[rr] && \cd^- A\ar@/_1pc/[ll]^{}\ar@/_-1pc/[ll]^{}\ar[rr] && \cd^- C.\ar@/_1pc/[ll]\ar@/_-1pc/[ll]}
\]
In particular, $\cd^-H$ is the triangle quotient of $\cd^-A$ by $\cd^-C$.

We claim that $P=eA$ and $Q=[0,Q';0]$ satisfy all the conditions of S.~K\"{o}nig's theorem (see subsection \ref{The mismatch}):
\emph{Condition 2.2:} $P$ is clearly an exceptional object of $\cd
A$ since it is a projective $A$-module.

\emph{Condition 2.3:}
\begin{itemize}
\item Let us check that $Q$ is exceptional. Since the canonical functor $\cd H\ra\cd A$ is fully faithful, we just have to check that $Q'$ is an exceptional object of $\cd H$, which is true because $Q'$ is an injective $H$-module.\item Since $H$ is hereditary, $Q'$ admits a projective resolution of length $1$. But the canonical functor $\Mod H\ra \Mod A$ preserves projective objects, and thus $Q$ admits a projective resolution of length $1$. This shows that $Q$ is isomorphic in $\cd A$ to a bounded complex of projective $A$-modules.\item To check that $\Hom_{A}(Q,?)$ preserves small coproducts of copies of $Q$one uses the fact that $\cd H\ra\cd A$ is fully faithful and applies Lemma\ref{preserves coproducts} with $Q'_H$.
\end{itemize}

\emph{Condition 2.4:} Since $P$ is a projective $A$-module, we only have to check that $(\cd A)(P,Q)=0$, but this is clear since\[(\cd A)(P,Q)\cong\Hom_{A}(P,Q)\cong Qe=0.
\]

\emph{Condition 2.5:} Let $M$ be a complex of $A$-modules such that $(\cd A)(P[i],M)=(\cd A)(Q[i],M)=0$ for each integer $i$. Consider the triangle\[\tau_{\cx}M\ra M\ra\tau^{\cy}M\ra(\tau_{\cx}M)[1]
\]
of $\cd A$. Since $(\cd A)(P,\tau^{\cy}M[i])=0$ for each integer $i$, then $(\cd A)(P,\tau_{\cx}M[i])=0$ for each integer $i$. Now, the fact that $P$ generates $\cx$ implies $\tau_{\cx}M=0$, that is to say, $M$ belongs to $\cy$. Therefore, we can assume that $M$ is the image of a complex $M'$ of $H$-modules by the canonical functor $\cd H\ra\cd A$. Then, since $H$ is right hereditary, for each integer $i$ we have\begin{align}
0=(\cd A)(Q[i],M)\cong(\cd H)(Q'[i],M')\cong(\cd H)(Q'[i],\prod_{n\in\Z}\H n(M')[-n])\cong \nonumber \\\cong \prod_{n\in\Z}(\cd H)(Q'[i],\H
n(M')[-n])\cong\prod_{n\in\Z}\Ext^{-n-i}_{H}(Q',\H n(M')) \nonumber
\end{align}
Finally, Lemma \ref{is a generator} (applied with $Q'_H$) tells us that
$M'$ is acyclic, that is to say, $M=0$ in $\cd A$.

According to S.~K\"{o}nig's theorem, $\cd^-A$ is a recollement as follows:\[\xymatrix{\cd^-(\End_{H}(Q'))\ar[rr] && \cd^- A\ar@/_1pc/[ll]^{}\ar@/_-1pc/[ll]^{}\ar[rr] && \cd^- C.\ar@/_1pc/[ll]\ar@/_-1pc/[ll]}
\]
In particular, $\cd^-(\End_{H}(Q'))$ is the triangle quotient of $\cd^-A$by $\cd^-C$. Therefore, $\cd^-(\End_{H}(Q'))$ is triangle equivalent to
$\cd^-H$. Let us fix a triangle equivalence
$F:\cd^-(\op{End}_{H}(Q'))\arr{\sim}\cd^-H$ and let us put
$F(\op{End}_{H}(Q'))=:T$. Since $H$ is right hereditary and $T$ is a
compact object of $\cd H$, we deduce that $T$ is isomorphic in $\cd H$ toa finite coproduct of stalk complexes $M_{i}[n_{i}]\ko 1\leq i\leq r\ko
n_{i}\in\Z$, for some $H$-modules $M_{i}$. This implies that each $M_{i}$is compact in $\cd H$. Therefore, each $M_{i}$ is finitely presented and
so (\cf \cite[Proposition I.12.1, Corollary I.11.5]{Stenstrom}) it is a
finitely generated projective $H$-module. Assume that $r>1$, and, withoutloss of generality, that $M_{i}\neq 0$ for each $1\leq i\leq r$ and that
$n_{i}\neq n_{j}$ for two different indexes $i$ and $j$. Since $T$ is
exceptional, there exists an isomorphism of algebras $\op{End}_{\cd
H}(T)\cong\bigoplus_{i=1}^{r}\op{End}_{\cd H}(M_{i})$ inducing a triangleequivalence $\cd H\simeq\bigoplus_{i=1}^{r}\cd(\op{End}_{\cd H}(M_{i}))$.This implies (\cf \cite[Example 1.7.15]{NicolasTesis}) that there exists anon-zero central idempotent $e$ of $H$ different from $1$, which
contradicts the fact that $H$ is a simple ring.
projective $H$-module. Of course, $T$ generates the triangulated category$\cd^-H$ and so it is also a generator of the abelian category $\Mod H$.
We have deduced that $T$ is a finitely generated projective generator of
$\Mod H$ and so $\op{End}_{H}(T)$ is Morita equivalent to $H$. In
particular, since $\op{End}_{H}(Q')\cong\op{End}_{\cd
H}(T)\cong\op{End}_{H}(T)$, we have that $\op{End}_{H}(Q')$ is Morita
equivalent to $H$. By the explicit description of Morita equivalences,
this is impossible because $H$ is countable and $\op{End}_{H}(Q')$ is not.
\begin{remark}
 S.~K\"{o}nig has pointed out to us that the construction of the functor $F$in \cite[Theorem 2.12]{Rickard1989} still yields a full embedding if $T$
is a bounded complex of (not necessarily finitely generated) projective
$A$-modules such that $(\cd A)(T,T^{(S)}[i])=0$ for every set $S$ and
every non-zero integer $i$ and that, as a consequence,  his  proof of
\cite[Theorem 1]{Konig1991} should still work assuming our hypothesis 2.3)of Theorem \ref{ourKonig}.
\end{remark}


\begin{thebibliography}{111}
\bibitem{AndersonFuller1992} F. W. Anderson and K. R. Fuller, \emph{Rings and Categories of Modules}, 2nd edition, Springer-Verlag Graduate Texts in Mathematics, \textbf{13}, Berlin, 1992.\bibitem{AuslanderReitenSmalo1995} M. Auslander, I. Reiten and S. O. Smalø, \emph{Representation Theory of Artin algebras},Cambridge St. Adv. Maths. \textbf{36}, Cambridge Univ. Press (1995).
\bibitem{BeilinsonBernsteinDeligne} A. A. Beilinson, J. Bernstein and P. Deligne, \emph{Faisceaux Pervers}, Ast\'{e}risque \textbf{100} (1982).\bibitem{BokstedtNeeman1993} M. B\"{o}kstedt and A. Neeman, \emph{Homotopy limits in triangulated categories}, Compositio Mathematica, \textbf{86} no. 2 (1993),  209--234.\bibitem{CisinskiNeeman2005} D.-C. Cisinski and A. Neeman, \emph{Additivity for derivator K-theory}, preprint (2005) available at Denis-Charles Cisinski's homepage.\bibitem{DinhGuilLopez} H. Q. Dinh, P. A. Guil Asensio and S. R. L\'{o}pez-Permouth, \emph{On the Goldie dimension of rings and modules}, J. Algebra {\bf 305} (2006), no. 2, 937--948.\bibitem{GeigleLenzing} W. Geigle and H. Lenzing, \emph{Perpendicular Categories with Applications to Representations and Sheaves}, Journal of Algebra \textbf{144} (1991), 273--343.\bibitem{GrusonJensen} L. Gruson and C. U. Jensen, \emph{Dimensions cohomologiques reli\'{e}es aux foncteurs $\op{lim}^{(i)}$}. Springer LNM {\bf 867}, 234--294 (1981).\bibitem{Hirschhorn2003} P. S. Hirschhorn, \emph{Model Categories and Their Localizations}, Math. Surveys and Monographs,\textbf{99}, American Mathematical Society, Providence, RI, 2003.
\bibitem{Keller1990} B. Keller, \emph{Chain complexes and stable categories}, Manus. Math. \textbf{67} (1990), 379--417.\bibitem{Keller1994a} B. Keller, \emph{Deriving DG categories}, Ann. Scient. Ec. Norm. Sup. (4) \textbf{27} (1994), no. 1, 63--102.\bibitem{Keller1998b} B. Keller, \emph{On the construction of triangle equivalences}, chapter of: \emph{Derived equivalences for group rings},Springer Lecture Notes in Mathematics \textbf{1685} (1998), 155--176, edited by Steffen K\"{o}nig and Alexander Zimmermann.\bibitem{Keller2006b} B. Keller, \emph{On differential graded categories}, in \emph{International Congress of Mathematicians.Vol. II}, pages 151--190. Eur. Math. Soc., Zürich, 2006.
\bibitem{KellerVossieck1987} B. Keller and D. Vossieck, \emph{Sous les cat\'{e}gories d\'{e}riv\'{e}es}, C. R. Acad. Sci. Paris S\'{e}r. I Math.\textbf{305} (1987), no. 6, 225--228.
\bibitem{KellerVossieck88b} B. Keller and D. Vossieck, \emph{Aisles in derived categories}, Bull. Soc. Math. Belg. S\'{e}r. A \textbf{40} (1988), no. 2, 239--253.\bibitem{Konig1991} S. K\"{o}nig, \emph{Tilting complexes, perpendicular categories and recollements of derived module categories of rings}, Journal ofPure and Applied Algebra \textbf{73} (1991), 211--232.
\bibitem{Krause2002} H. Krause, \emph{A Brown representability theorem via coherent functors}, Topology \textbf{41} (2002), 853--861.\bibitem{Maltsiniotis2001} G. Maltsiniotis, \emph{Introduction \`{a} la th\'{e}orie des d\'{e}rivateurs (d'apr\`{e}s Grothendieck)}, preprint (2001) available at Georges Maltsiniotis' homepage.\bibitem{Maltsiniotis2005} G. Maltsiniotis, \emph{La K-th\'{e}orie d'un d\'{e}rivateur triangul\'{e} (suivi d'un appendice par B. Keller)}, ``Categories in Algebra, Geometry and Mathematical Physics'', Contemp. Math. \textbf{431} (2007), 341--368.\bibitem{Milnor1962} J. Milnor, \emph{On Axiomatic Homology Theory}, Pacific J. Math. \textbf{12} (1962), 337--341.\bibitem{Neeman2001} A. Neeman, \emph{Triangulated categories}, Annals of Mathematics studies. Princeton University Press, \textbf{148}, (2001).\bibitem{Nicolas2007} P. Nicol\'{a}s, \emph{The bar derived category of a curved dg algebra}, arXiv:math/0702449v1 [math.RT], to appear in Journal of Pure and Applied Algebra.\bibitem{NicolasSaorin2007c} P. Nicol\'{a}s and M. Saor\'{i}n, \emph{Parametrizing recollement data}, submitted.\bibitem{NicolasTesis} P. Nicol\'{a}s, \emph{On torsion torsionfree triples}, Ph. D. Thesis, Universidad de Murcia, 2007.\bibitem{Porta2007} M. Porta, \emph{The Popescu-Gabriel theorem for triangulated categories}, arXiv: math.KT/0706.4458v1.\bibitem{Rickard1989} J. Rickard, \emph{Morita theory for Derived Categories}, J. London Math. Soc. \textbf{39} (1989), 436--456.\bibitem{Souto2004} M. J. Souto Salorio, \emph{On the cogeneration of t-structures}, Arch. Math. \textbf{83} (2004), 113--122.\bibitem{Stenstrom} B. Stenstr\"om, \emph{Rings of quotients}, Grundlehren der math. Wissensch., \textbf{217}, Springer-Verlag, 1975.\end{thebibliography}
\end{document}